\documentclass[francais]{smfart}
\makeindex
\usepackage{epsfig,amsmath,a4,amssymb,amsfonts,amsbsy,xy,latexsym,epsf,pstricks,multicol,accents}
\xyoption{all}
\usepackage[french]{babel}
\theoremstyle{plain}

\let\set\mathbb
\def\<<{\leavevmode
  \raise0.28ex\hbox{$\scriptscriptstyle\langle\!\langle$}\nobreak
  \hskip -.6pt plus.3pt minus.2pt\,}
\def\>>{\,\nobreak\hskip -.6pt plus.3pt minus.2pt
  \raise0.28ex\hbox{$\scriptscriptstyle\rangle\!\rangle$}}

\newcommand{\indexentry}[2]{{\item #1,  page #2}}

\def\bB{{\|}}

\def\SL2Z{{{\rm SL}_2(\ZZ)}}

\def\ba{{\bf a}}
\def\bff{{\bf f}}
\def\bb{{\bf b}}
\def\bh{{\bf h}}

\def\bk{{\bf k}}

\def\bzero{{{\bf 0}_g}}
\def\bzerod{{{\bf 0}_{g_2}}}

\def\bun{{{\bf 1}_g}}
\def\bund{{{\bf 1}_{g_2}}}
\def\bdeux{{{\bf 2}_g}}

\def\bn{{\bf n}}

\def\bl{{\bf    l}}

\def\bm{{\bf m}}
\def\bc{{\bf c}}

\def\br{{\bf r}}

\def\bx{{\bf x}}
\def\bq{{\bf q}}
\def\bz{{\bf z}}
\def\by{{\bf y}}

\def\AA{{\set A}}

\def\CC{{\set C}}
\def\Div{{\rm Div}}

\def\FF{{\set F}}

\def\Fq{{\FF _q}}

\def\NN{{\set N}}
\def\PP{{\set P}}
\def\PP{{\set P}}
\def\PU{{{\set P}^1}}
\def\Pic{\mathop{\rm{Pic}}\nolimits }
\def\QQ{{\set Q}}

\def\RR{{\set R}}

\def\ZZ{{\set Z}}

\def\cA{{\mathcal A}}
\def\cB{{\mathcal B}}

\def\cD{{\mathcal D}}

\def\cH{{\mathcal H}}

\def\cJ{{\mathcal J}}
\def\cK{{\mathcal K}}
\def\cL{{\mathcal L}}

\def\cM{{\mathcal M}}
\def\cN{{\mathcal N}}
\def\cO{{\mathcal O}}
\def\cP{{\mathcal P}}
\def\cQ{{\mathcal Q}}
\def\cR{{\mathcal R}}
\def\cS{{\mathcal S}}

\def\sc{\scriptscriptstyle }

\newtheorem{lemme}{Lemme}

\newtheorem{theoreme}{Th{\'e}or{\`e}me}

\newtheorem{definition}{D{\'e}finition}

\author[Jean-Marc Couveignes]{Jean-Marc  Couveignes}
\address{Groupe de Recherche en Informatique et Math\'ematiques  du Mirail\\
Universit\'e de Toulouse II, Le Mirail}
\email{couveig@univ-tlse2.fr}
\urladdr{http://www.univ-tlse2.fr/grimm/couveignes}
\title{Jacobiens, jacobiennes  et stabilit{\'e} num{\'e}rique}
\alttitle{Jacobians and numerical stability}

\begin{document}

\frontmatter

\begin{abstract}
On {\'e}tudie la complexit{\'e} et la stabilit{\'e}  des calculs dans la jacobienne des courbes de grand genre sur
le corps des complexes avec une attention particuli{\`e}re aux courbes modulaires. 
\end{abstract}

\begin{altabstract}
This paper is concerned with  the complexity and stability of
arithmetic operations  in the jacobian variety
of curves over the field of complex numbers, as the genus grows to
infinity. We focus on modular curves.
\end{altabstract}

\subjclass{11F11, 11F25, 11F30, 11Y16, 11Y35, 65E05, 65Y20, 68Q15}

\keywords{jacobienne, approximation, stabilité, formes modulaires,
  complexité algorithmique, machine de Turing, temps polynomial déterministe}
\altkeywords{jacobian variety, complex approximation, stability,
  modular forms, algorithmic complexity, Turing machine, deterministic
polynomial time}

\bibliographystyle{smfplain}
\maketitle

\mainmatter
\tableofcontents

\section{Introduction}\label{section:introduction}

Il est traditionnel de calculer  dans le groupe des points de la jacobienne d'une courbe alg{\'e}brique projective lisse et g{\'e}om{\'e}triquement irr{\'e}ductible
$X$
de genre $g$
en repr{\'e}sentant tout {\'e}l{\'e}ment de ce groupe par un diviseur effectif  de degr{\'e} $g$, une fois choisi
un tel diviseur $O$ comme origine. La somme de deux diviseurs 
 $P-O$ et $Q-O$ est {\it r{\'e}duite} par le calcul
de l'espace lin{\'e}aire associ{\'e} au diviseur $P+Q-O$ suivi de la localisation
 des z{\'e}ros d'une fonction non nulle
de cet espace. 

Comme l'application de Jacobi,

$$S^gX\rightarrow J_X$$
\noindent 
de la puissance sym{\'e}trique $g$-i{\`e}me $S^gX$ \index{$S^gX$, la puissance sym{\'e}trique $g$-i{\`e}me de $X$} de $X$ dans sa jacobienne $J_X$\index{$J_X$, la jacobienne de $X$}, n'est pas un
isomorphisme, la repr{\'e}sentation n'est pas unique.

Si le corps de base est un corps fini $\Fq$, les op{\'e}rations
arithm{\'e}tiques y
 sont exactes et rapides.

On consid{\`e}re ici   le cas o{\`u} le corps de base est le corps $\CC$ des complexes.
On se donne un mod{\`e}le analytique naturel et une mesure sur $X(\CC)$.
On s'int{\'e}resse {\`a}  la complexit{\'e}  des algorithmes utilis{\'e}s pour ajouter
et r{\'e}duire des diviseurs. Le cadre est celui des machines de Turing 
classiques. En effet, on peut avoir
 en vue des applications arithm{\'e}tiques comme le calcul de
nombres
de points, ou de coefficients de formes modulaires 
 et les calculs en nombres complexes ne sont alors qu'une {\'e}tape dans la recherche d'une quantit{\'e} discr{\`e}te.
Le projet de Bas Edixhoven pour r{\'e}pondre {\`a} une question de Ren{\'e} Schoof \cite{edix2, edix1, couv} 
se prête à cette approche.

Bien s{\^u}r, les machines de Turing ordinaires ne manipulent  pas les nombres r{\'e}els ni complexes mais
plut{\^o}t des nombres rationnels, d{\'e}cimaux ou binaires.
Cependant, on peut voir un nombre r{\'e}el $\alpha$  comme un oracle qui, pour tout entier
positif $k$, retourne une valeur binaire ou d{\'e}cimale approch{\'e}e 
de $\alpha$ {\`a} $\exp(-k)$ pr{\`e}s. Si une machine de Turing doit r{\'e}soudre un probl{\`e}me
dont les entr{\'e}es sont des nombres r{\'e}els, elle re{\c c}oit un oracle
pour chacun de ces r{\'e}els. Si la machine de Turing calcule un nombre r{\'e}el, on lui donne
en entr{\'e}e la pr{\'e}cision absolue $k$ requise et elle retourne une valeur approch{\'e}e {\`a} $\exp(-k)$ pr{\`e}s 
du r{\'e}sultat. On dit que la machine est polynomiale si elle r{\'e}pond en temps polynomial en la taille
des donn{\'e}es et $k$. On note que la recherche des racines complexes d'un
polyn{\^o}me unitaire {\`a} coefficients complexes se fait en temps d{\'e}terministe polynomial
gr{\^a}ce {\`a}  la m{\'e}thode de quadrichotomie de Weyl par exemple.  On veut dire par l{\`a}
qu'une valeur approch{\'e}e {\`a}  $\exp(-k)$ pr{\`e}s  de chaque racine peut {\^e}tre calcul{\'e}e en temps
polynomial
en le degr{\'e} du polyn{\^o}me, la taille des coefficients (logarithme du maximum des modules
des coefficients) et la pr{\'e}cision absolue $k$ requise.

On veut savoir si la complexit{\'e} asymptotique des op{\'e}rations arithm{\'e}tiques dans
la jacobienne est polynomiale en le genre de la courbe. La premi{\`e}re difficult{\'e} est de donner un sens
pr{\'e}cis {\`a} cette assertion. Plut{\^o}t que de rester dans le vague, on formule et on
{\'e}tudie ces
questions
dans le cas important et repr{\'e}sentatif des courbes modulaires $X_0(p)$ lorsque $p$ est un entier premier qui
tend
vers l'infini. L'algorithmique de ces courbes est riche et largement explor{\'e}e. On trouve
dans  \cite{cremona, delaunay, elkies} des algorithmes  pour l'{\'e}tude
homologique des courbes modulaires et des m{\'e}thodes analytiques exp{\'e}rimentales
motiv{\'e}es par la v{\'e}rification de conjectures arithm{\'e}tiques et la recherche de points rationnels.

La section \ref{section:modules}  d{\'e}crit le mod{\`e}le analytique standard de ces  courbes ainsi
que ses propri{\'e}t{\'e}s algorithmiques. On y rappelle d'abord  les r{\'e}sultats de Manin, Shokurov, Cremona 
et Merel concernant le calcul des p{\'e}riodes, et on en donne une expression quantifi{\'e}e du point de vue
de la complexit{\'e} algorithmique et de la stabilit{\'e} num{\'e}rique. Cette derni{\`e}re est assur{\'e}e en dernier ressort
par des minorations du volume des p{\'e}riodes et du d{\'e}terminant
jacobien de l'application d'int{\'e}gration de Jacobi.
Ces minorations reposent elles m{\^e}mes sur des consid{\'e}rations d'int{\'e}gralit{\'e} des coefficients des formes
modulaires primitives, propres et normalis{\'e}es.

On pr{\'e}sente dans la section \ref{section:arithmetique}  des algorithmes pour les op{\'e}rations {\'e}l{\'e}mentai\-res dans la jacobienne
$J_0(p)$ et pour la r{\'e}solution effective du probl{\`e}me inverse de Jacobi. La complexit{\'e} et la stabilit{\'e}
de ces algorithmes sont {\'e}tudi{\'e}es avec les outils de la section \ref{section:modules} puis estim{\'e}es dans les
th{\'e}or{\`e}mes \ref{theoreme:arithmetiquedebase} et
\ref{theoreme:inverse}. On obtient des algorithmes déterministes
polynomiaux en $p$. Le caractère déterministe de ces algorithmes
s'explique en dernier lieu par la connexité  du tore analytique complexe $J_0(p)(\CC)$.

Tous les lemmes et d{\'e}finitions concernant  la localisation et la stabilit{\'e} des z{\'e}ros de fonctions analytiques
sont pr\'esent{\'e}s dans l'appendice  \ref{section:analyse} qui est 
ind{\'e}pen\-dant mais doit {\^e}tre au moins parcouru avant
de lire les sections \ref{section:modules} et \ref{section:arithmetique}.

Les m{\'e}thodes, les {\'e}nonc{\'e}s  et les démonstrations que nous donnons pour les courbes $X_0(p)$ s'{\'e}tendent 
sans peine
au cas de $X_1(p)$.
Pour les courbes modulaires de niveau compos{\'e},  il faut une majoration 
des coefficients  des  d{\'e}veloppe\-ments de Fourier en toutes les pointes 
ainsi qu'un algorithme pour les calculer.

On trouvera un index {\`a} la fin de cet article.

{\bf Convention  importante : } le symbole $\cO$ d{\'e}signe partout une constante absolue positive
et effective,  chaque fois 
diff{\'e}rente. La pr{\'e}sence de ce symbole dans une formule ou un {\'e}nonc{\'e}
signifie que cette formule ou cet {\'e}nonc{\'e} sont vrais si,  pour chaque occurence, ce symbole est remplac{\'e}e par une 
constante positive effective
bien choisie.

\section{Courbes modulaires $X_0(p)$}\label{section:modules}

Cette section rappelle, pr{\'e}cise et compl{\`e}te quelques r{\'e}sultats m{\'e}triques et algorithmiques concernant les
courbes modulaires $X_0(p)$. On supposera que $p$ est premier et que  le genre $g$ de $X_0(p)$
est au moins $2$.

Le paragraphe \ref{subsection:lesdeuxdisques} introduit quelques notations et un recouvrement 
non injectif de $X_0(p)$
par deux disques analytiques centr{\'e}s en chacune des deux pointes.

Les propri{\'e}t{\'e}s {\'e}l{\'e}mentaires des formes primitives, propres et normalis{\'e}es sont rappel{\'e}es dans le paragraphe \ref{subsection:differentielles}
et celles de l'homologie dans le paragraphe \ref{subsection:homologie}. Le calcul
des p{\'e}riodes est abord{\'e} dans le paragraphe \ref{subsection:periodes}.  Ces trois paragraphes r{\'e}sument le travail
de Manin, Shokurov, Cremona et Merel sur cette question. 

Le paragraphe \ref{subsection:torecomplexe} {\'e}tablit une minoration du volume du r{\'e}seau des p{\'e}riodes. Une formule d'int{\'e}gration
sur les surfaces de Riemann relie ce volume au produit des normes de Petersson des formes primitives, propres et normalis{\'e}es, ces derni{\`e}res
{\'e}tant faciles {\`a} minorer parce que le d{\'e}veloppement de
Fourier commence par $(1+O(q))dq$ o\`u $q=\exp(2i\pi\tau)$ est le
param\`etre de Tate associ{\'e} \`a un
 $\tau$ du demi-plan de Poincar{\'e}.

Le paragraphe \ref{subsection:jacobi} {\'e}tablit des majorations simples mais n{\'e}cessaires des int{\'e}grales de Jacobi et 
d{\'e}finit l'{\it instabilit{\'e} } d'un diviseur effectif de degr{\'e} $g$.
Le paragraphe \ref{subsection:jacobiensetwronskiens} construit un diviseur d'instabilit{\'e} assez petite.  Cela revient
{\`a} trouver $g$ points dans le mod{\`e}le canonique de $X_0(p)$ qui ne soient proches d'aucun  hyperplan.
Autrement dit,  le jacobien en ces $g$ points n'est pas trop petit. On prend le parti (maladroit en pratique
mais simple en th{\'e}orie) de chercher les $g$ points dans le voisinage de la pointe {\`a} l'infini. Le terme principal
du d{\'e}veloppement du jacobien y est le wronskien. On le minore gr{\^a}ce {\`a} l'int{\'e}gralit{\'e} des coefficients
de son d{\'e}veloppement de Fourier.

Le paragraphe \ref{subsection:stabilite} {\'e}tudie la stabilit{\'e} de l'application inverse de Jacobi. Cel{\`a} se r{\'e}duit {\`a} majorer
la diff{\'e}rence entre cette application et sa lin{\'e}aris{\'e}e.

La connaissance d'un $g$-uplet de points de faible instabilit{\'e}, donn{\'e} au  paragraphe \ref{subsection:jacobiensetwronskiens},
permet de construire au paragraphe \ref{subsection:reseaux} des sous-ensembles finis de taille modeste et bien distribu{\'e}s
 dans  le tore complexe.
Comme les {\'e}l{\'e}ments de ces ensembles sont images par l'application de Jacobi de diviseurs connus, il sont des
auxiliaires pr{\'e}cieux pour la r{\'e}solution approch{\'e}e du
probl{\`e}me inverse de Jacobi. Ils permettent de discrétiser ce problème.

\subsection{Un mod{\`e}le analytique}\label{subsection:lesdeuxdisques}

Soit $p$ un nombre premier et $X= X_0(p)$
la courbe modulaire de niveau $p$ associ{\'e}e au sous groupe de congruence   $\Gamma = \Gamma_0(p)$ de $\SL2Z$.
On note $\cH$\index{$\cH$, le demi-plan de Poincar{\'e}} le demi-plan de Poincar{\'e} et $\cH^*=\cH\cup\QQ\cup \{\infty\}$.
\index{$\cH^*$, le demi-plan de Poincar{\'e} avec les pointes}  La surface de Riemann compacte
quotient 
$\Gamma  \backslash \cH^*$ est   $X(\CC)$. Son genre est $g=\frac{p+1-3\nu_2-4\nu_3}{12}$ 
avec $\nu_2=1+\left(  \frac{-1}{p}\right)$\index{$\nu_2$, le nombre de points elliptiques d'ordre $2$}
 et $\nu_3= 1+\left(  \frac{-3}{p} \right)$. \index{$\nu_3$, le nombre de points elliptiques d'ordre $3$}
Il y a $\nu_2$ points elliptiques d'ordre $2$ et on note $\cP_2$\index{$\cP_2$, le diviseur des points
elliptiques d'ordre $2$}  le diviseur somme de ces
points. De m{\^e}me il a $\nu_3$ points elliptiques d'ordre $3$ et on note $\cP_3$\index{$\cP_3$, le diviseur des points
elliptiques d'ordre $3$} le diviseur somme de ces
points. 
Voir \cite[Propositions 1.40 et 1.43]{Shi}.  Le genre de $X$   est compris entre $\frac{p-13}{12}$ et  $\frac{p+1}{12}$.
 Le quotient $\Gamma \backslash \cH$  
est un ouvert de Zariski de $X$ not{\'e} $Y=Y_0(p)$. 
On note que  la largeur de la pointe $\infty$ est $1$ et la largeur de la pointe $0$ est
$p$. Pour $\tau \in \cH$  on pose $q=q(\tau)=q_\infty(\tau)= \exp(2i\pi\tau)$\index{$q$, $q_0$, $q_\infty$, $P_0$, $P_\infty$} et 
$w(\tau)=-\frac{1}{p\tau}$ et $q'=q'(\tau)=q_0(\tau)=q(w(\tau))= \exp(\frac{-2i\pi}{p\tau})$.
On note $P=P_\infty=P(\tau)=P(q)$ le point de $Y$ associ{\'e} {\`a} $\tau$ et $P'=P_0=P'(\tau)=P'(q)=P(w(\tau))=W(P)=P(q')$
o{\`u} $W$\index{$W(P)$, l'involution d'Atkin-Lehner appliqu{\'e}e au point $P$} est l'involution d'Atkin-Lehner. On a le diagramme

\begin{equation*}\xymatrix{
Y \ar@{->}[rr]^W&&Y \\
D-\{0\} \ar@{->}[urr]^{P_0}   \ar@{->}[u]^{P_\infty}&&D-\{0\} \ar@{->}[u]_{P_\infty}\\
\cH \ar@{->}[u]^q \ar@{->}[rr]^w  \ar@{->}[urr]^{q'} &&\cH   \ar@{->}[u]_q}
\label{diagramme:gros}
\end{equation*}

\'Etant donn{\'e}s deux  r{\'e}els $R_\infty$ et $R_0$  plus petits que $1$ on peut se demander si
l'union de l'image par $P_\infty$ du disque ouvert  $D(0,R_\infty)$\index{$D(x,r)$, le disque de centre $x$ et de rayon
$r$} et de l'image par 
$P_0$ de $D(0,R_0)$ recouvre $X(\CC)$.

On pose $S= \left( \begin{array}{cc} 0 & -1 \\ 1 & 0   \end{array} \right)$  et $T= \left( 
\begin{array}{cc} 1 & 1 \\ 0 & 1   \end{array} \right)$ de sorte que 
 $S\tau=-1/\tau$ et $T\tau=\tau+1$. Soit  $\cR$ le domaine fondamental usuel de $\SL2Z$, d{\'e}limit{\'e} 
par le cercle de centre $0$ et de rayon $1$ et par les droites d'abscisses $-1/2$ et $1/2$.
Alors un domaine fondamental pour $\Gamma$ est constitu{\'e}  de l'union de $\cR$ et des
$ST^k\cR$ pour $k$ entier de  $0$ {\`a} $p-1$. Ces  derniers sont contenus dans 
l'image par $S$ de l'ensemble des $\tau =a+ib$ avec $b\ge \frac{\sqrt{3}}{2}$. Donc leur image par $w$
est constitu{\'e}e de complexes dont la partie imaginaire est au moins $\frac{\sqrt{3}}{2p}$. Si on 
choisit $R_0>\exp(-\frac{\pi\sqrt{3}}{p})$ alors l'image de $D(0,R_0)$ par $P_0$
recouvre les $ST^k\cR$ pour $k$ de $0$ {\`a} $p-1$.

Comme $\cR$ est contenu dans le demi plan des parties imaginaire au moins {\'e}gales {\`a} $\sqrt{3}/2$
on prend $R_\infty>\exp(-\pi\sqrt{3})$ et l'image de $D(0,R_\infty)$ par $P_\infty$
recouvre $\cR$.

On pose donc
$R_\infty=0.005$\index{$R_\infty = 0.005$ } et  $R_0=1-\frac{1}{p}$\index{$R_0 = 1-\frac{1}{p}$ }

On a donc recouvert $X(\CC)$ par l'image de deux disques analytiques
complexes $D_\infty=D(0,R_\infty)$ et $D_0=D(0,R_0)$.

\subsection{Diff{\'e}rentielles}\label{subsection:differentielles}

On peut maintenant calculer des espaces de formes diff{\'e}rentielles sur  $X$. On fixe
donc un entier $d\ge 1$. 
\`A toute  forme modulaire  parabolique $f$ de  poids $2d$  sur $\Gamma$
on associe la diff{\'e}rentielle $\omega = (2i\pi)^d f(d\tau)^d$ de degr{\'e} $d$.
D'apr{\`e}s \cite[Proposition 2.16]{Shi} on a 

$$\Div (\omega)=\Div(f)-d(0) -d(\infty)-\frac{d}{2}\cP_2 
-\frac{2d}{3}\cP_3.$$

On pose donc $\Delta_d= (d-1)(0) + (d-1)(\infty)+\lfloor \frac{d}{2}\rfloor \cP_2 
+ \lfloor \frac{2d}{3}\rfloor \cP_3$\index{$\Delta_d$, le diviseur $(d-1)(0) + (d-1)(\infty)+\lfloor \frac{d}{2}\rfloor \cP_2 
+ \lfloor \frac{2d}{3}\rfloor \cP_3$} et on cherche 
 une base $\cD_d$ de l'espace 
$\cH^{d}(  \Delta_d)$\index{$\cH^{d}(  \Delta_d)$, l'espace des formes 
diff{\'e}rentielles de degr{\'e} $d$ et de diviseur $\ge -\Delta_d$} des formes 
diff{\'e}rentielles de degr{\'e} $d$ et de diviseur $\ge -\Delta_d$.

On prend pour $\cD_d$\index{$\cD_d$, l'ensemble des $\omega =(2i\pi)^d f(q)(d\tau)^d=\frac{f(q)}{ q^d}(dq)^d$ 
o{\`u} $f(q)$ est une forme modulaire parabolique primitive, propre et
normalis{\'e}e  sur $\Gamma=\Gamma_0(p)$ et de poids $2d$}
 l'ensemble des $\omega =(2i\pi)^d f(q)(d\tau)^d=\frac{f(q)}{ q^d}(dq)^d$ 
o{\`u} $f(q)$ est une forme modulaire parabolique primitive\footnote{Cela signifie qu'elle ne provient pas d'une forme de niveau plus petit. Se dit en Anglais
  {\it newform}. Cette condition est vide ici puisque le niveau $p$
  est premier.}, propre\footnote{Autrement dit, $f$ est vecteur propre
  des op{\'e}rateurs de Hecke.} et  normalis{\'e}e\footnote{Son
  d{\'e}veloppement de Fourier commence par $q$.} sur $\Gamma=\Gamma_0(p)$
et de poids $2d$.
Si $f$ est une telle  forme  elle admet un d{\'e}veloppement 
$f=\sum_{k\ge 1}a_k q_{\infty}^k$ avec $a_1=1$ et pour tout entier $k\ge 1$ on montre que
le coefficient 
$a_k$ est un entier alg{\'e}brique major{\'e} en module par $k^{d+2}$.
Il suffit de le montrer pour $k=\ell ^n$ une puissance d'un premier $\ell$.
D'apr{\`e}s  le th{\'e}or{\`e}me de Deligne   on a $|a_\ell|\le 2 \ell^{d-\frac{1}{2}}$ 
 et d'apr{\`e}s \cite[Theorem 3]{AtkinLehner} 

$$|a_{\ell^{n+2}}|\le |a_{\ell}a_{\ell^{n+1}}|+\ell^{2d-1}|a_{\ell^{n}}|$$
\noindent 
donc $|a_{\ell^{n}}|\le u_n\ell^\frac{n(2d-1)}{2}$  o{\`u} $u_n$ est la suite r{\'e}currente $u_0=1$, $u_1=2$ et
$u_{n+2}=2u_{n+1}+u_n$. Donc $u_n=\frac{(1+\sqrt 2)^{n+1}-(1-\sqrt 2)^{n+1}}{2\sqrt 2}$ et
$|u_n|\le 4^n\le \ell ^{2n}$ donc $|a_{\ell^n}|\le \ell^{2n}\ell^\frac{n(2d-1)}{2}$.

Le d{\'e}veloppement de $\omega$ en $q_\infty$ est donc le d{\'e}veloppement standard, donn{\'e} par les valeurs
propres des op{\'e}rateurs de Hecke.  On peut calculer les coefficients $a_k$
comme valeurs propres des op{\'e}rateurs de Hecke agissant  sur les symboles de Manin-Shokurov suivant 
\cite{cremona, merel,frey}. Les plongements complexes des valeurs propres
peuvent alors {\^e}tre approch{\'e}s en temps polynomial en $p$ et
la pr{\'e}cision absolue requise, par un algorithme de recherche de racines de polyn{\^o}mes comme celui de Weyl.

Le d{\'e}veloppement de $\omega$ en $q_0$ est le tir{\'e} en arri{\`e}re de $\omega$ par l'application
$P_0 : D-\{0\} \rightarrow Y$. Comme $P_0$ est la compos{\'e}e de $P_\infty$ et de $W$, le d{\'e}veloppement
de $\omega$ en $q_0$ est le d{\'e}veloppement de $W(\omega)$ en $q_\infty$. Mais 
$\omega$ est vecteur propre de $W$ de valeur propre $\pm1$. On a donc la m{\^e}me majoration
pour les coefficients du d{\'e}veloppement de $\omega$ en $q_0$. 

\begin{lemme}[Manin, Shokurov, Cremona, Merel]

Il existe un algorithme  qui pour tous $p$ premier, $d\ge 1$ et $r \ge 1$ calcule les plongements complexes
des $r$ premiers coefficients de toutes les formes modulaire primitives, propres et normalis{\'e}es de niveau $p$ et poids $2d$
en temps polynomial en $p$, $d$,  $r$ et la pr{\'e}cision absolue requise.
\end{lemme}

\subsection{L'homologie de la courbe}\label{subsection:homologie}

La th{\'e}orie de Manin \cite{manin,merel,cremona} {\'e}tablit
que l'homologie relative $H_1(X,ptes,\ZZ)$ est engendr{\'e}e par les symboles modulaires. Un symbole
est not{\'e} indiff{\'e}rement $(c:d)=\left( \begin{array}{cc} a & b \\ c & d   \end{array}   \right)=\{\frac{b}{d}, \frac{a}{c}  \}$.
  L'ensemble des symboles
est $\PP=\PU (\ZZ/p \ZZ)$. 
On rappelle que  $S= \left( \begin{array}{cc} 0 & -1 \\ 1 & 0   \end{array} \right)$  et $T= \left( 
\begin{array}{cc} 1 & 1 \\ 0 & 1   \end{array} \right)$. Pour tout $\gamma \in \SL2Z$ on note $(\gamma )$ le symbole
$\{\gamma (0), \gamma (\infty) \}$.

On note ${\bf B}$ le 
sous-$\ZZ$-module   de ${\ZZ}^{\PP}$  
engendr{\'e} par les  $(c:d)+(c:d)S=(c:d)+(-d:c)$ et $(c:d)+(c:d)TS+(c:d)(TS)^2=(c:d)+(c+d:-c)+(d:-c-d)$ o{\`u} 
$(c:d)$ parcourt $\PP$.

On note ${\bf Z}$ le sous $\ZZ$-module libre et satur{\'e}
de $\ZZ^{\PP}$ engendr{\'e} par les $(c:1)$  pour $c\not=0$ et par $(\infty)=(0:1)+(1:0)$. C'est
le module des symboles {\`a}  bord nul. La base form{\'e}e des   $(c:1)$ pour $c\not = 0$ et de $(\infty)$
permet d'identifier {\bf Z} au r{\'e}seau  $\ZZ^{p}$ de $\RR^{p}$
muni de la forme bilin{\'e}aire canonique (de matrice identit{\'e} dans 
cette base). 

On note que ${\bf B} \subset {\bf Z} \subset \RR^p$.  Comme le quotient ${\bf Z}/{\bf B}=H_1(X,\ZZ)$ est
sans-torsion, le sous module ${\bf B}$ est satur{\'e} dans ${\bf Z}$. On identifie $H_1(X,\ZZ)$ {\`a} la projection orthogonale
de ${\bf Z}$ sur le $\RR$-espace vectoriel de $\RR^p$ orthogonal au sous-espace vectoriel  $\RR{\bf B}$ engendr{\'e}
par ${\bf B}.$

Comme ${\bf B}$ est engendr{\'e} par des vecteurs de norme $\le \sqrt 3$ et qu'il a pour dimension
$p-2g$,  son volume $V$ est un entier positif major{\'e} par $3^\frac{p-2g}{2}$. 
Mais d'apr{\`e}s  \cite[Proposition I.2.9, Proposition I.3.5]{martinet}
$H_1(X,\ZZ)$ est inclus dans $\frac{1}{V^2}\ZZ^p$ et  son  volume  est {\'e}gal {\`a} $1/V$.

D'apr{\`e}s
l'in{\'e}galit{\'e} d'Hermite  \cite[Th{\'e}or{\`e}me II.2.1]{martinet}, le r{\'e}seau  $V^2H_1(X,\ZZ)\subset \ZZ^p$ admet une base  constitu{\'e}e 
de vecteurs de norme\footnote{Noter que dans le livre de Martinet, la norme 
d'un vecteur est
d\'efinie comme la valeur de la forme quadratique
en ce vecteur. Ici, la norme est la racine carr\'ee de la forme
quadratique.}  major{\'e}e par $\left(\frac{4}{3}\right)^\frac{g(g-1)}{2}
3^\frac{p-2g}{2}$. L'algorithme LLL \cite[Theorem 2.6.2]{cohen} 
produit en temps polynomial  en $p$ une base de $V^2H_1(X,\ZZ)$ form{\'e}e de vecteurs entiers 
de norme $\le 2^{\frac{g(g-1)}{4}}3^\frac{p-2g}{2}\le 3^{p^2}$. Donc ces vecteurs sont des combinaisons de symboles
avec des coefficients major{\'e}s par cette m{\^e}me borne. On peut faire beaucoup mieux en y regardant de plus
pr{\`e}s.
Au total, on obtient une base $\cB$ de $H_1(X,\ZZ)\subset \ZZ^p$  dont les vecteurs sont des combinaisons
de symboles {\`a} coefficients dans $\frac{1}{V^2}\ZZ$ dont les num{\'e}rateurs sont des entiers born{\'e}s en valeur absolue
par $3^{p^2}$.

\subsection{Les p{\'e}riodes}\label{subsection:periodes}

On observe  que le diviseur $\Delta_1$ est nul. Donc les formes primitives, propres et  normalis{\'e}es
de poids $2$ donnent une base $\cD_1$ de  l'espace  des formes diff{\'e}rentielles holomorphes. 

Le r{\'e}seau $\cL$ des p{\'e}riodes est construit en int{\'e}grant chaque forme de $\cD_1$
le long des symboles comme font Tingley dans sa th{\`e}se  et Cremona dans son livre
\cite[Proposition 2.10.1]{cremona}. Si $g \in \Gamma_0(p)$
 avec 
$g = \left( \begin{array}{cc} a & b \\ pc & d   \end{array}   \right)$ et $c>0$
on pose $y_0=\frac{1}{pc}$ et $x_1=-dy_0$ et $x_2=ay_0$. Alors
pour toute forme $f=\sum_{k\ge 1}a_k q_{\infty}^k$ de poids $2$ on a 

\begin{equation}\label{equation:integration}
\int_{0}^{g(0)} f(q) \frac{dq}{q}= 
\sum_{n\ge 1} \frac{a_n}{n}\exp(-2\pi n y_0) (\exp(2\pi inx_2 ) - \exp( 2\pi i n x_1)  ).
\end{equation}

Cette quantit{\'e} est {\'e}valu{\'e}e en temps polynomial en $p$, $c$ et la pr{\'e}cision absolue requise et elle est major{\'e}e en
module par un polyn{\^o}me en $p$ et $c$.

Si $(c:1)=\left( \begin{array}{cc} 1 & 0 \\ c & 1   \end{array}   \right)=\{0, \frac{1}{c}  \}$
est un symbole diff{\'e}rent de $(0:1)$ et $(1:0)$, on peut supposer
que $1\le c < p$. 
 Afin d'utiliser la formule
d'int{\'e}gration (\ref{equation:integration}), on choisit deux entiers
$u$ et $v$ tels que $uc-vp=1$ et on note que 
la matrice  $g=   \left( \begin{array}{cc} u & 1 \\ pv & c   \end{array}   \right) \in \Gamma_0(p)$ v{\'e}rifie $g(0)=1/c$.
On peut choisir $0\le u, v<p$. 
Donc $\int_{(c:1)}f(q) \frac{dq}{q}$ est {\'e}valu{\'e}e en temps
polynomial en $p$ et la pr{\'e}cision absolue requise,
 et elle est 
major{\'e}e par un polyn{\^o}me en $p$.

 Comme les {\'e}l{\'e}ments de la base $\cB$ de $H_1(X,\ZZ)$ construite ci-dessus sont des
combinaisons lin{\'e}aires des symboles $(c:1)$ et $(\infty)\in {\bf B}$  avec des coefficients  dans $\frac{1}{V^2}\ZZ$ 
{\`a} num{\'e}rateurs major{\'e}s par $3^{p^2}$ en valeur absolue,
on peut calculer les p{\'e}riodes $\int_\gamma \omega$ pour $\gamma \in \cB$ et $\omega \in \cD_1$ en temps polynomial
en $p$ et la pr{\'e}cision absolue requise et ces p{\'e}riodes sont major{\'e}es en module par $\exp(p^{\cO})$.

\subsection{Le tore complexe}\label{subsection:torecomplexe}

L'espace $H=\cH^1\oplus \bar \cH^1$\index{$H=\cH^1\oplus \bar \cH^1$, l'espace des  formes harmoniques}  des formes harmoniques admet une forme bilin{\'e}aire
d{\'e}finie par int{\'e}gration. Si $\omega_1 = u_1 dq + \overline{v_1 dq  }$ et  $\omega_2 = u_2 dq + \overline{v_2 dq  }$ 
on pose $<\omega_1, \omega_2>=\int_X \omega_1 \wedge  \omega_2$.

L'int{\'e}gration d{\'e}finit  aussi un accouplement $\int : H\times
H_1(X,\CC)\rightarrow \CC$  qui à  $(\omega, \gamma)$ associe
la p{\'e}riode $\int_\gamma \omega$. 
Il en r{\'e}sulte un isomorphisme entre $H$ et le dual de $H_1(X,\CC)$. Mais l'accouplement d'intersection
induit un isomorphisme entre $H_1(X,\CC)$ et son dual.
On en d{\'e}duit un isomorphisme $\iota$ entre $H$ et $H_1(X,\CC)$ qui {\`a} tout $\omega$ associe l'unique $\gamma = \iota (\omega)$
tel que $\int_g \omega=\gamma \cdot  g$ pour tout $g\in H_1(X,\CC)$. D'apr{\`e}s
\cite[Proposition III.2.3.]{farkas}, cet isomorphisme est une isom{\'e}trie :

$$<\omega_1, \omega_2>=\int_X \omega_1 \wedge  \omega_2= \iota (\omega_1) \cdot  \iota (\omega_2).$$

On d{\'e}finit l'op{\'e}rateur
${}^* : H \rightarrow H$ par ${}^*(u dq + \overline{v dq  })=(-iu dq + i\overline{v dq  })$.
On d{\'e}finit un produit hermitien sur $H$ par 

$$(\omega_1, \omega_2)=\int_X \omega_1 \wedge {}^*\bar \omega_2=i\int_X (u_1\bar u_2+\bar v_1 v_2) dq\wedge \overline{ dq}.$$

C'est le produit de Petersson.
Les op{\'e}rateurs de Hecke sont autoadjoints pour ce produit hermitien. Donc deux formes distinctes dans $\cD_1$ sont 
orthogonales. Il reste {\`a} {\'e}valuer $(\omega, \omega)$ pour chaque {\'e}l{\'e}ment $\omega$ de la base $\cD_1$.

On note $\tilde \cB$ la base de $H_1(X,\ZZ)$ duale {\`a} gauche de $\cB$ pour la forme d'intersection, c'est-{\`a}-dire que pour
tout $\gamma \in \cB$ il existe un unique  $\tilde \gamma \in \tilde \cB$ tel que $\tilde \gamma \cdot \gamma =1$
et si $\gamma \not = \gamma'$ on a $\tilde \gamma \cdot \gamma' =0$.
Ainsi  $\iota (\omega)= \sum_{\gamma \in \cB} \tilde \gamma \int_\gamma\omega$.

Donc

$$(\omega, \omega)=<\omega, {}^*\bar \omega>=\iota (\omega) \cdot  \iota ({}^*\bar \omega)= \iota (\omega) \cdot  \iota (i\bar \omega)=
i\sum_{\gamma, \gamma ' \in \cB} \tilde \gamma \cdot  \tilde \gamma' \int_\gamma \omega \int_{\gamma'}\bar \omega .$$

Notons $P$ la matrice $2g\times g$ des p{\'e}riodes holomorphes

$$P=\left( \int_\gamma \omega  \right)_{\gamma \in \cB,\, \omega \in \cD_1.}$$

Soit $M=(P | \bar P)$ la matrice $2g\times 2g$ des p{\'e}riodes harmoniques. Soit $M^*=(i\bar P | -i  P)$.
Soit $\cQ=\left( \gamma \cdot  \gamma'  \right)_{\gamma, \gamma' \in \cB}$ la matrice de la forme d'intersection dans la base $\cB$.
Soit $\tilde \cQ={}^t\cQ^{-1}=\left( \tilde \gamma \cdot  \tilde \gamma'  \right)_{\gamma, \gamma' \in \cB}$ la matrice de la forme d'intersection dans la base $\tilde \cB$.

La matrice du produit scalaire de Petersson dans la base $\cD_1 \cup
 \bar \cD_1$ de $H$ n'est autre   que
 ${}^tM\tilde\cQ M^*$.
Comme le d{\'e}terminant de $\cQ$ est $1$ et comme  le volume du r{\'e}seau $\cL$ des p{\'e}riodes est le 
module du  d{\'e}terminant de $M$ divis{\'e} par $2^g$, il vient que
ce volume est le produit des $\frac{1}{2}(\omega, \omega)$ pour $\omega$ dans  la base $\cD_1$.

On peut minorer chacun des $(\omega, \omega)$ en notant que 

$$\omega = \frac{f(q)}{q} dq=(1+\sum_{k\ge 2}a_kq^{k-1})dq$$
\noindent 
avec $|a_k|\le k^3$ de sorte que $\int_X  \omega \wedge {}^*\bar \omega =\int_X |1+\sum_{k\ge 2}a_kq^{k-1}|^2idq\wedge\overline{dq}$.
On observe que le disque form{\'e} des $q$ de module inf{\'e}rieur {\`a}   $\exp(-2\pi)$ est contenu dans un domaine fondamental de $X$
donc $(\omega, \omega)\ge 2\pi r^2\min_{|q|\le r} |1+\sum_{k\ge 2}a_kq^{k-1}|^2$ pour tout $r\le \exp(-2\pi)$.

Or
pour $|q|\le r$ on a $|\sum_{k\ge 1 }a_{k+1}q^k|\le r\sum_{k\ge 0} (k+2)^3r^k\le \frac{8r}{(1-r)^4}\le 0.016$ si $r=\exp(-2\pi)$.
Ainsi $(\omega, \omega)\ge 2\pi \exp(-4\pi)(1-0.016)^2\ge 2\cdot 10^{-5}$ de sorte que le volume du r{\'e}seau $\cL$ des p{\'e}riodes est
au moins $10^{-5g}$.

Comme les p{\'e}riodes $\int_\gamma \omega$ sont major{\'e}es par $\exp(p^{\cO})$ on en d{\'e}duit que le volume de tout sous-r{\'e}seau du r{\'e}seau
des p{\'e}riodes est minor{\'e} par $\exp(-p^{\cO})$. De sorte que si l'on conna{\^\i}t les p{\'e}riodes avec une pr{\'e}cision absolue polynomiale en
$p$,
on les conna{\^\i}t aussi avec une bonne pr{\'e}cision relative. En particulier, si on conna{\^\i}t un point de  $\CC^{\cD_1}$
par ses coordonn{\'e}es complexes alors on conna{\^\i}t le point du tore $\CC^{\cD_1}/\cL$ puisque $\cL$ n'est pas trop petit ni trop
aplati.

\begin{lemme}[Volume et complexit{\'e} du r{\'e}seau des p{\'e}riodes]\label{lemme:lereseaudesperiodes}
Si $X_0(p)$ est de genre $g\ge 1$ on note $\cD_1$ la base de 
$\cH^1(X_0(p))$ constitu{\'e}e des formes primitives, propres et  normalis{\'e}es
et on identifie le dual de $\cH^1(X_0(p))$ {\`a} $\CC^{\cD_1}$.
On appelle r{\'e}seau des p{\'e}riodes le r{\'e}seau de $\CC^{\cD_1}$
form{\'e} des p{\'e}riodes de $X_0(p)$. Ce r{\'e}seau
est de volume $\ge 10^{-5g}$. Tous les sous-r{\'e}seaux non nuls du r{\'e}seau des p{\'e}riodes ont un volume $\ge \exp( -p^{c_1})$
o{\`u} $c_1$ est une constante positive effective.
Le r{\'e}seau des p{\'e}riodes admet une base constitu{\'e}e de vecteurs de norme $\le \exp(p^{c_2})$ o{\`u} $c_2$
est une constante positive effective. Une telle base  peut {\^e}tre calcul{\'e}e
en temps polynomial en $p$ et la pr{\'e}cision absolue requise.
\end{lemme}

\subsection{L'application d'int{\'e}gration de Jacobi}\label{subsection:jacobi}

On note $\mu_0 : D(0,R_0)\rightarrow \CC^{\cD_1}$ l'application d'int{\'e}gration de Jacobi

$$\mu_0 : q_0\mapsto \left(\int^{q_0}\omega\right)_{\omega\in \cD_1.}$$

 On ne pr{\'e}cise pas l'origine $o$
de l'int{\'e}grale. Cette application est bien d{\'e}finie {\`a} une constante
additive pr{\`e}s. On d{\'e}finit de m{\^e}me $\mu_\infty :  D(0,R_\infty)\rightarrow \CC^{\cD_1}$ en veillant {\`a} choisir la m{\^e}me origine $o$

$$\mu_\infty : q_\infty\mapsto \left(\int^{q_\infty}\omega\right)_{\omega\in \cD_1.}$$

Ces int{\'e}grales se calculent par int{\'e}gration terme {\`a} terme de la s{\'e}rie associ{\'e}e {\`a} la forme diff{\'e}rentielle.

Le disque $D(0,1)\subset \CC$ est muni de la distance usuelle associ\'ee
\`a la norme sur $\CC$ d\'efinie par le   module $z\mapsto |z|$.

L'espace $\CC^{\cD_1}$ peut \^etre muni des normes $L^2$ ou $L^\infty$
dont la d\'efinition est rappel\'ee au paragraphe
\ref{subsection:ordredegrandeur} de l'appendice.

L'application d'int{\'e}gration de Jacobi est Lipschitzienne dans  le sens suivant.

Soit $P_1$ un point de coordonn{\'e}es $\tau_1$, $q_1$ et $q_1'$ tel que $q_1=\exp(2i\pi \tau_1)$ est 
dans $D(0,R_\infty)=D(0,0.005)$. Soit $P_2$ de coordonn{\'e}es  $\tau_2$, $q_2$ tel que $q_2$ est 
proche de $q_1$, en ce sens que  $q_2\in D(0,0.01)$. Pour tout $\omega \in \cD_1$ on a
$|\int_{q_1}^{q_2}\omega|\le |q_2-q_1|\max_{|q|\le 0.01} (|\sum_{k\ge 1}a_kq^{k-1}|)$.
Or pour $|q|\le 0.01$ on a 
$|\sum_{k\ge 1}a_kq^{k-1}|\le \sum_{k\ge 0} (k+1)^3 10^{-2k} \le \frac{6}{0.99^4} \le 7$.
Donc 

$$|\mu_\infty(P_2)-\mu_\infty(P_1)|_\infty \le 7 |q_2-q_1| \mbox{ et } |\mu_\infty(P_2)-\mu_\infty(P_1)|_2\le 7\sqrt g |q_2-q_1|.$$

Soit maintenant $P_1$ un point de coordonn{\'e}es $\tau_1$, $q_1$ et $q_1'$ tel que $q_1'=
\exp(\frac{-2i\pi}{p\tau_1})$ est 
dans $D(0,R_0)=D(0,1-\frac{1}{p})$. Soit $P_2$ de coordonn{\'e}es  $\tau_2$, $q_2'$ tel que $q_2'$ est 
proche de $q_1'$, en ce sens que  $q_2'\in D(0,1-\frac{1}{2p})$. Pour tout $\omega \in \cD_1$ on a
$|\int_{q_1'}^{q_2'}\omega|\le |q_2'-q_1'|\max_{|q|\le 1-\frac{1}{2p}} (|\sum_{k\ge 1}a_kq^{k-1}|)$.
Or pour $|q|\le 1- \frac{1}{2p}$ on a 
$|\sum_{k\ge 1}a_kq^{k-1}|\le \sum_{k\ge 0} (k+1)^3 (1-\frac{1}{2p})^{k} \le {96p^4}$.
Donc 

$$|\mu_0(P_2)-\mu_0(P_1)|_\infty \le 96p^4 |q_2'-q_1'| \mbox{ et } |\mu_0(P_2)-\mu_0(P_1)|_2\le 96p^4\sqrt g |q_2'-q_1'|.$$

Ainsi la perte de pr{\'e}cision occasionn{\'e}e par l'application de Jacobi est $O(\log p)$.

\begin{lemme}[Majoration des int{\'e}grales de Jacobi]
Pour tout   premier $p$ on pose $R_\infty=0.005$ et $R_0=1-\frac{1}{p}$ et on 
recouvre  $X_0(p)$ par les  deux disques analytiques $D_\infty$ et $D_0$ 
centr{\'e}s en chacune des deux pointes $\infty$ et $0$ et de rayons respectifs $R_\infty$ et $R_0$. 
Donc  $D_\infty= D(0,R_\infty)=\{q_\infty , |q_\infty|\le 0.005 \}$ et 
$D_0= D(0,R_0)=\{q_0 , |q_0|\le 1-\frac{1}{p}\}$. 

Sur chacun de ces deux disques, l'int{\'e}gration de Jacobi d{\'e}finit une application
{\`a} valeur dans le dual de l'espace des formes holomorphes $\cH^1(X_0(p))$. 
On munit   $\cH^1(X_0(p))$ de la base $\cD_1$ constitu{\'e}e 
des formes primitives, propres et  normalis{\'e}es et on note abusivement $\CC^{\cD_1}$ son dual, que l'on munit
de  la norme    $L^\infty$  associ{\'e}e  {\`a} la base canonique (duale de $\cD_1$). 

L'application d'int{\'e}gration de Jacobi  est  alors Lipschitzienne sur chacun des deux 
disques (et même sur leurs voisinages $D(0,0.01)$ et
$D(0,1-\frac{1}{2p})$) et son coefficient de dilatation y  est 
major{\'e} par $7$ sur le premier  et par $96p^4$ sur le second.
\end{lemme}

Soit $\epsilon = (\epsilon_k)_{1\le k \le g}\in \{0,\infty \}^g$. 
On note $D_\epsilon $\index{$D_\epsilon $, le produit 
$D_\epsilon =  D(0,R_{\epsilon_1})\times
\cdots \times  D(0,R_{\epsilon_g})$ pour un $\epsilon = (\epsilon_k)_{1\le k \le g}\in \{0,\infty \}^g$} le produit 

$$D_\epsilon =  D(0,R_{\epsilon_1})\times
\cdots \times  D(0,R_{\epsilon_g}).$$

Les  $D_\epsilon$ recouvrent le produit $X(\CC)^g$.
L'application  produit \index{$\mu_\epsilon$ pour $\epsilon =  (\epsilon_k)_{1\le k \le g}\in \{0,\infty \}^g$}

$$\mu_\epsilon = \mu_{\epsilon_1}\times \cdots\times
\mu_{\epsilon_g} : D_\epsilon
\rightarrow \CC^{\cD_1}$$
\noindent  associe au $g$-uplet $(q_1,\ldots,q_g)$ la 
somme des $\mu_{\epsilon_k}(q_{\epsilon_k,k})$ pour $1\le k \le g$.  

On note $(z_\omega)_\omega$ la base duale
de  la base canonique de $\CC^{\cD_1}$. On note $\cJ_{\cD_1,\epsilon}$ le  d{\'e}terminant jacobien de 
l'application $\mu_\epsilon$ en $(q_{\epsilon_1,1},\ldots,q_{\epsilon_g,g})$ \index{$\cJ_{\cD_1,\epsilon}$, le  d{\'e}terminant jacobien de 
l'application $\mu_\epsilon$ en $(q_{\epsilon_1,1},\ldots,q_{\epsilon_g,g})$} :

$$\cJ_{\cD_1,\epsilon}(q_{\epsilon_1,1},\ldots,q_{\epsilon_g,g})=
\frac{\bigwedge z_\omega}{\bigwedge dq_{\epsilon_k,k}}(q_{\epsilon_1,1},\ldots,q_{\epsilon_g,g})=
 \left| \frac{\omega}{dq_{\epsilon_k}}(q_{\epsilon_k,k})\right|_{\omega,k.}$$

Ce d{\'e}terminant est
une s{\'e}rie enti{\`e}re de $g$ variables dont les coefficients
se majorent ais{\'e}ment {\`a} partir d'une majoration des coefficients
des formes $\omega$. Mais il  peut parfaitement s'annuler.
L'{\it instabilit{\'e}} d'un $g$-uplet dans $D_\epsilon$  peut  se d{\'e}finir
comme l'oppos{\'e}  du logarithme du module de ce d{\'e}terminant jacobien.

\subsection{Jacobiens et wronskiens}\label{subsection:jacobiensetwronskiens}

L'espace $\cH^1$ des diff{\'e}rentielles holomorphes est muni de la base $\cD_1$ et de la norme $L^\infty$ associ{\'e}e. Si $F=
\sum_{\omega \in \cD_1} f_\omega \omega$ cette norme est not{\'e}e $|F|=|f|=\max_\omega|f_\omega|.$
Puisqu'on dispose de deux disques analytiques sur $X(\CC)$ centr{\'e}s en chacune des deux pointes, il est 
naturel d'introduire une norme $\cS_\infty(F)=\max_{q_\infty\in \bar D(0,R_\infty)}\left| \frac{F}{dq_\infty}\right|$ et
de m{\^e}me $\cS_0(F)=\max_{q_0\in \bar
  D(0,R_0)}\left|\frac{F}{dq_0}\right|$. On introduit aussi les
variantes  $\hat \cS_\infty(F)=\max_{q_\infty\in \bar D(0, \frac{1}{2} )}\left| \frac{F}{dq_\infty}\right|$ et
de m{\^e}me $\hat \cS_0(F)=\max_{q_0\in \bar
  D(0,\frac{1}{2})}\left|\frac{F}{dq_0}\right|$.
\index{$\cS_\infty(F)$,  $\cS_0(F)$,  $|F|$, $\hat \cS_\infty(F)$, et 
  $\hat \cS_0(F)$,
les cinq normes de la forme $F$}

Puisque toutes les normes sont {\'e}quivalentes, les quotients
$\cS_0(F)/|F|$, $\cS_\infty (F)/|F|$, $\hat \cS_0(F)/|F|$, $\hat\cS_\infty (F)/|F|$ pour $F\not=0$ sont major{\'e}s
et minor{\'e}s par des bornes ind{\'e}pendantes de $F$. On veut montrer que le logarithme de ces bornes
est polynomial en $p$. C'est {\'e}vident pour la borne sup{\'e}rieure parce que les formes de la base 
$\cD_1$ s'{\'e}crivent $\omega=f(q_\infty)dq_\infty=\pm f(q_0)dq_0$
avec $f$ d'ordre de grandeur\footnote{La d\'efinition de l'ordre de
  grandeur d'une s\'erie est donn\'ee au paragraphe \ref{subsection:ordredegrandeur}.} $(1,3)$. Il suffit d'appliquer
le lemme \ref{lemme:majorationdureste} de majoration du reste.

Pour controler la borne inf{\'e}rieure, il suffit de trouver $g$ points $q_1$, \ldots, $q_g$ de 
$D(0,R_\infty)$ tels que le module
du jacobien $\cJ_{\cD_1,(\infty, \ldots, \infty )}(q_1,\ldots,q_g)$ ait un logarithme  born{\'e} inf{\'e}rieurement
par un polyn{\^o}me $-g^{\cO}$  en le genre $g$.

Une strat{\'e}gie possible est de chercher d'abord un $q$ tel que le wronskien de $\cD_1$ en $q$ ne soit pas trop
petit et de chercher ensuite des $q_1$,\ldots , $q_g$ dans un voisinage de ce $q$.

On utilise le classique 

\begin{lemme}[wronskien et jacobien]\label{lemme:WJ}
Soient $g\ge 2$ un entier naturel et $f_1(q)$, $f_2(q)$, \ldots, $f_g(q)$ des s{\'e}ries de Laurent  {\`a} coefficients
complexes. On appelle wronskien associ{\'e} {\`a} $\bff = (f_1, \ldots, f_g)$  le d{\'e}terminant \index{$W_\bff $, le wronskien
associ{\'e} {\`a} la famille $\bff$}

$$W_\bff(q)=\left|  \begin{array}{ccc} f_1(q) & \dots & f_g(q)\\ f'_1(q) & \dots & f'_g(q)\\ \vdots & &\vdots \\
f^{(g-1)}_1(q) & \dots & f^{(g-1)}_g(q)\end{array} \right|.$$

On se donne $g$ ind{\'e}termin{\'e}es $q_1$, $q_2$, \ldots, $q_g$ et on appelle jacobien associ{\'e} {\`a} $\bff$ le d{\'e}terminant

$$\cJ_\bff = \left|  \begin{array}{ccc} f_1(q_1) & \dots & f_g(q_1)\\ f_1(q_2) & \dots & f_g(q_2)\\ \vdots & &\vdots \\
f_1(q_{g}) & \dots & f_g(q_g)\end{array} \right|.$$\index{$\cJ_\bff $, le jacobien
associ{\'e} {\`a} la famille $\bff$}

On note $D=\prod_{k<l}(q_l-q_k)$ le discriminant r{\'e}duit.

Si les $f_k$ sont des s{\'e}ries enti{\`e}res alors
 le jacobien $\cJ_\bff$ est dans l'anneau $\CC[[q_1,\ldots,q_g]]$ et il
est divisible par le discriminant  r{\'e}duit $D$ dans cet anneau. Le quotient $\cJ_\bff/D$ est alors congru {\`a}
$W_\bff(0)$ modulo l'id{\'e}al maximal de $\CC[[q_1,\ldots,q_g]]$.
\end{lemme}

La démonstration  peut se faire par r{\'e}currence sur $g$. \hfill $\Box$

\begin{lemme}[Majoration du wronskien]\label{lemme:ordrewronskien}
Si pour $1\le l \le g$ les $f_l(q)$ sont des s{\'e}ries enti{\`e}res d'ordre de grandeur 
$(A,k)$ avec $A\ge 1$ et $k\ge 1$ alors le wronskien $W_\bff(q)$ est
une s{\'e}rie enti{\`e}re d'ordre de grandeur 

$$( g! A^g 2^{\frac{g(g-1)(g+3k-2)}{6}}, gk-1+\frac{g(g+1)}{2}).$$
\end{lemme}

Cela r{\'e}sulte du lemme \ref{lemme:derivprod}. \hfill $\Box$

\begin{lemme}[Majoration du jacobien]\label{lemme:ordrejacobien}
Si pour $1\le l\le g$ les $f_l(q)$ sont des s{\'e}ries enti{\`e}res d'ordre de grandeur 
$(A,k)$ avec $A\ge 1$ et $k\ge 1$ alors le jacobien $\cJ_\bff(q_1,\ldots,q_g)$ est
une s{\'e}rie enti{\`e}re des $g$ variables $q_1$, \ldots, $q_g$ d'ordre de grandeur 

$$( g! A^g, (k,k,\ldots,k)).$$
\end{lemme}

Selon le  lemme \ref{lemme:WJ}, le wronskien nous renseigne  sur
 le terme dominant du jacobien au voisinage 
de $q_1=q_2=\cdots=q_g=0$.

\begin{definition}[Wronskien et Jacobien des formes diff{\'e}rentielles]
Pour  chaque forme diff{\'e}rentielle primitive, propre et  normalis{\'e}e $\omega \in \cD_1$  on rappelle que  $\omega=\frac{f(q)}{q}dq$ o{\`u} $f$
est une forme modulaire primitive, propre et  normalis{\'e}e de poids $2$. On forme
la famille $\bff$ des $g$ s{\'e}ries enti{\`e}res $\frac{f(q)}{q}$ ainsi obtenues. Pour
tout $q\in D(0,1)$, le wronskien
de $\bff$ en $q$ est not{\'e} $W_\infty (q)$ \index{ $W_\infty (q)$, le wronskien associ{\'e}
aux d{\'e}veloppements de Fourier {\`a} l'infini des formes diff{\'e}rentielles}  et le jacobien correspondant en $\bq=(q_1, \ldots, q_g)$  est not{\'e} $\cJ_\infty(\bq)$.
\index{$\cJ_{\infty}(\bq)$, le jacobien associ{\'e} {\`a} la famille
des d{\'e}veloppements {\`a} l'infini des formes diff{\'e}rentielles}
\end{definition}

 Posons $m=\frac{g(g+1)}{2}$. 
Le produit $W_\infty (q)(dq)^m$ 
est une forme holomorphe de degr{\'e} $m$ sur $X(\CC)$.
 Donc  elle  a au plus $2(g-1)m$ z{\'e}ros en comptant
les multiplicit{\'e}s. En particulier, la pointe {\`a} l'infini $q_\infty=0$ est un
z{\'e}ro de multiplicit{\'e} $v \le 2(g-1)m$ de $W_\infty (q)$.

La s{\'e}rie enti{\`e}re $W_\infty (q)^2$ a des coefficients entiers rationnels et sa valuation $2v$
est au plus $4(g-1)m=2g(g^2-1)$. On suppose $g\ge 2$.
 Le lemme \ref{lemme:ordrewronskien} de majoration
du wronskien montre que l'ordre de grandeur de $W_\infty (q)^2$ est $(\exp(g^{\cO}),g^{\cO})$. Le lemme \ref{lemme:majorationdureste} 
de majoration du reste permet d'{\'e}crire
$W_\infty (q)^2=wq^{2v}+R_{2v+1}(q)$ avec $w$ entier naturel non nul et $R_{2v+1}$
major{\'e} en module par $|q|^{2v+1}(2v+2)^{g^{\cO}}$
pour $q\in D(0,R_\infty)$. Donc  $W_\infty (q)^2$ est minor{\'e} en module 
par $\frac{1}{2}q^{2g(g^2-1)}$ si $|q|\le (2v+2)^{-g^{\cO}}.$

\begin{lemme}[Minoration du wronskien]\label{lemme:minorationduwronskien}
Il existe une constante effective positive $c_6$ telle que pour tout premier $p$
tel que  le genre  $g$ de $X_0(p)$ est au moins $2$, il existe un $q\in D(0,10^{-10})$
tel que $\log \left|W_\infty (q)\right| \ge -g^{c_6}.$ Un tel $q$ se calcule en temps polynomial
en $p$ et la pr{\'e}cision absolue requise.
\end{lemme}

Soit $q\in D(0, R_\infty)$ tel que $W_\infty(q)\not = 0$ et posons $\bq=(q, \ldots, q)$ et $\bx = (x_1, \ldots, x_g)$. 
Le jacobien $\cJ_\infty(\bq+\bx\star(\bun -\bB \bq \bB ))$ est\footnote{La
  notation $\star$ est introduite au paragraphe \ref{subsection:ordredegrandeur}.} une s{\'e}rie enti{\`e}re en les $g$ variables $x_1$, \ldots, $x_g$. 
C'est en fait le jacobien
des $(f_k(q+x(1-|q|))/(q+x(1-|q|)))_{1\le k \le g}$ vues comme fonctions de $x$. Le lemme \ref{lemme:WJ} appliqu{\'e} {\`a} ces derni{\`e}res donne le terme
principal de cette s{\'e}rie en $\bx=\bzero$ :

$$\cJ_\infty(\bq+\bx \star (\bun -\bB \bq \bB ) )=W_\infty(\bq)(1-|q|)^{\frac{g(g-1)}{2}}\prod_{k<l}(x_l-x_k) +R_{\frac{g(g-1)}{2}+1}(\bx).$$

Le jacobien $\cJ_\infty(\bq)$  est  une s{\'e}rie en $\bq$ 
d'ordre de grandeur $(g!,(3,3,\ldots,3))$ d'apr{\`e}s le lemme \ref{lemme:ordrejacobien}.
Si $q\in D(0,R_\infty)$, la s{\'e}rie recentr{\'e}e
$\cJ_\infty(\bq+\bx\star(\bun - \bB \bq \bB ))$ est une s{\'e}rie en $\bx$ d'ordre
de grandeur $(\exp(g^{\cO}), (4,4, \ldots,4))$ d'apr{\`e}s le lemme \ref{lemme:plusrecentrage} de recentrage.
Pour $\bx $ dans $P(\bzero ,R_\infty)$, le lemme \ref{lemme:majorationdureste} de majoration du reste donne alors pour
$g\ge 2$ 

$$\left| R_{\frac{g(g-1)}{2}+1}(\bx)\right| \le \exp(g^{\cO})|\bx|_\infty^{\frac{g(g-1)}{2}+1}.$$

On pose $s=|\bx|_\infty$ et on suppose que $\bx=(\frac{s}{g},\frac{2s}{g}, \ldots, \frac{(g-1)s}{g}, s)$ et $s\le R_\infty$.
Alors

\begin{eqnarray*}
\left|W_\infty(\bq)(1-|q|)^{\frac{g(g-1)}{2}}\prod_{k<l}(x_l-x_k)\right| &\ge& \left| W_\infty(\bq)\right|
\left(\frac{s(1-|q|)}{g}\right)^{\frac{g(g-1)}{2}}\\
&\ge & \left| W_\infty(\bq)\right| \left(\frac{0.995s}{g}\right)^{\frac{g(g-1)}{2}}.
\end{eqnarray*}

On choisissant un $\bq$ donn{\'e} par le lemme \ref{lemme:minorationduwronskien} et en prenant
$s=\exp(-g^{\cO})$ on prouve le 

\begin{lemme}[Minoration du jacobien]\label{lemme:minorationdujacobien}
Il existe une constante effective $c_7$
telle que pour tout premier $p$ tel que  le genre  $g$ de
 $X_0(p)$ est au moins $2$, il existe  $\br= (r_1, \ldots, r_g) \in D(0,R_\infty)^g$
tel que $\log \left|\cJ_\infty(\br)\right| \ge -g^{c_7}.$ Un tel $\br$ se calcule en temps polynomial
en $p$ et la pr{\'e}cision absolue requise.
Les  cinq normes sur $\cH^1$ d{\'e}finies pour $F=\sum_{\omega \in \cD_1} f_\omega \omega$ par $|F|=\max_\omega|f_\omega|$,
  $\cS_\infty(F)=\max_{q_\infty\in \bar
  D(0,R_\infty)}\left|\frac{F}{dq_\infty}\right|$, 
 $\cS_0(F)=\max_{q_0\in \bar D(0,R_0)}\left|\frac{F}{dq_0}\right|$,
  $\hat \cS_\infty(F)=\max_{q_\infty\in \bar
  D(0,\frac{1}{2})}\left|\frac{F}{dq_\infty}\right|$, et
 $\hat \cS_0(F)=\max_{q_0\in \bar D(0,\frac{1}{2})}\left|\frac{F}{dq_0}\right|$
 sont dans des rapports $\exp(p^{\cO})$,  ou, si l'on pr{\'e}f{\`e}re, il existe une constante positive  $c_3$  telle que pour tout $p$ et
pour tout $F\not = 0$  les 
diff{\'e}rences entre $\log |F|$, $\log \cS_0(F)$,  $\log \cS_\infty
(F)$, $\log \hat \cS_0(F)$ et   $\log \hat \cS_\infty (F)$, sont  born{\'e}es en valeur absolue
par $p^{c_3}$.
\end{lemme}

\subsection{Jacobiens et wronskiens quadratiques}

Des r{\'e}sultats semblables sont vrais pour le jacobien et le wronskien  associ{\'e}s {\`a} la
base $\cD_2$ de l'espace de  formes diff{\'e}rentielles quadratiques $\cH^2(\Delta_2)$.

\begin{definition}[Jacobien et wronskien des formes  quadratiques]
Pour   chaque forme diff{\'e}rentielle quadratique primitive, propre et  normalis{\'e}e $\phi \in  \cD_2$ on rappelle 
que  $\phi=\frac{h(q)}{q}\frac{(dq)^2}{q}$ o{\`u} $h$
est une forme modulaire primitive, propre et  normalis{\'e}e de poids $4$. On forme
la famille $\bh$ des $g_2 = 3g-1+\nu_2+\nu_3$\index{$g_2 = 3g-1+\nu_2+\nu_3$, la dimension
de $\cH^2(\Delta_2)$}  s{\'e}ries enti{\`e}res $\frac{h(q)}{q}$ ainsi obtenues. Pour
tout $q\in D(0,1)$, le wronskien
de $\bh$ en $q$ est not{\'e} $W_{2,\infty}(q)$\index{$W_{2,\infty}(q)$, le wronskien associ{\'e} {\`a} la famille
des d{\'e}veloppements {\`a} l'infini des formes quadratiques}  et le jacobien correspondant en $\bq=(q_1, \ldots, q_{g_2})$  est not{\'e} 
$\cJ_{2,\infty}(\bq)$.\index{$\cJ_{2,\infty}(\bq)$, le jacobien  associ{\'e} {\`a} la famille
des d{\'e}veloppements {\`a} l'infini des formes quadratiques}
\end{definition}

 Posons $m_2=\frac{g_2(g_2+3)}{2}$. \index{$m_2=\frac{g_2(g_2+3)}{2}$, le degr{\'e} du wronskien quadratique}
Le produit $W_{2,\infty}(q)q^{-g_2}(dq)^{m_2}$ 
est une forme  de degr{\'e} $m_2$ sur $X(\CC)$, holomorphe en dehors de $\Delta_2$. Plus pr{\'e}cis{\'e}ment
elle appartient {\`a} $\cH^{m_2}(g_2\Delta_2)$.
 Donc elle  a au plus $2(g-1)m_2+g_2(2+\nu_2+\nu_3)$ z{\'e}ros en comptant
les multiplicit{\'e}s. En particulier, la pointe {\`a} l'infini $q_\infty=0$ est un
z{\'e}ro de multiplicit{\'e} $v_2 \le 2(g-1)m_2+g_2(2+\nu_2+\nu_3)$ de $W_{2,\infty}(q)$. Donc $v_2 \le 9g(g+1)^2$.

La s{\'e}rie enti{\`e}re $W_{2,\infty }(q)^2$ a des coefficients entiers rationnels et sa valuation $2v_2$
est au plus $18g(g+1)^2$. On suppose $g\ge 2$.
 Commes les $\frac{h(q)}{q}$ sont d'ordre de grandeur $(1,4)$, le lemme \ref{lemme:ordrewronskien} de majoration
du wronskien montre que l'ordre de grandeur de $W_{2,\infty} (q)^2$ est $(\exp(g^\cO),g^\cO)$. Le lemme \ref{lemme:majorationdureste} 
de majoration du reste permet d'{\'e}crire
$W_{2,\infty } (q)^2=wq^{2v_2}+R_{2v_2+1}(q)$ avec $w$ entier naturel non nul et $R_{2v_2+1}$
major{\'e} en module par $|q|^{2v_2+1}(2v_2+2)^{g^\cO}$
pour $q\in D(0,R_\infty)$. Donc  $W_{2,\infty } (q)^2$ est minor{\'e} en module 
par $\frac{1}{2}q^{18g(g+1)^2}$ si $-\log |q|\ge g^\cO$ et $g\ge 2$.

\begin{lemme}[Minoration du wronskien quadratique]\label{lemme:minorationduwronskien2}
Il existe une cons\-tan\-te effective $c_8$ telle que pour tout premier $p$
tel que  le genre  $g$ de $X_0(p)$ est au moins $2$, il existe un $q\in D(0,10^{-10})$
tel que $\log \left|W_{2,\infty }(q)\right| \ge -g^{c_8}.$ Un tel $q$ se calcule en temps polynomial
en $p$ et la pr{\'e}cision absolue requise.
\end{lemme}

On obtient de m{\^e}me l'analogue du lemme \ref{lemme:minorationdujacobien} pour le jacobien des formes
quadratiques.

Soit $q\in D(0, R_\infty)$ tel que $W_{2,\infty}(q)\not = 0$ et posons $\bq=(q, \ldots, q)$ et $\bx = (x_1, \ldots, x_{g_2})$. 
Le jacobien $\cJ_{2,\infty}(\bq+\bx\star(\bund -\bB \bq \bB ))$ est une s{\'e}rie enti{\`e}re en les $g_2$ variables $x_1$, \ldots, $x_{g_2}$. 
C'est en fait le jacobien
des $(h_k(q+x(1-|q|))/(q+x(1-|q|)))_{1\le k \le g_2}$ vues comme fonctions de $x$. Le lemme \ref{lemme:WJ} appliqu{\'e} {\`a} ces derni{\`e}res donne le terme
principal de cette s{\'e}rie en $\bx=\bzerod$ :

$$\cJ_{2,\infty}(\bq+\bx)=W_{2,\infty}(\bq)(1-|q|)^{\frac{g_2(g_2-1)}{2}}\prod_{k<l}(x_l-x_k) +R_{\frac{g_2(g_2-1)}{2}+1}(\bx).$$

Le jacobien $\cJ_{2,\infty }(\bq)$  est  une s{\'e}rie en $\bq$ 
d'ordre de grandeur $(g_2!,(4,4,\ldots,4))$ d'apr{\`e}s le lemme \ref{lemme:ordrejacobien}.
Si $q\in D(0,R_\infty)$, la s{\'e}rie recentr{\'e}e
$\cJ_{2,\infty}(\bq+\bx \star (\bund - \bB \bq \bB ))$ est une s{\'e}rie en $\bx$ d'ordre
de grandeur $(\exp(g^{\cO}), (5,5, \ldots,5))$ d'apr{\`e}s le lemme \ref{lemme:plusrecentrage} de recentrage.
Pour $\bx $ dans $P(\bzero ,R_\infty)$, le lemme \ref{lemme:majorationdureste} de majoration du reste donne alors pour
$g\ge 2$ 

$$\left| R_{\frac{g_2(g_2-1)}{2}+1}(\bx)\right| \le \exp(g^\cO)|\bx|_\infty^{\frac{g_2(g_2-1)}{2}+1}.$$

On pose $s=|\bx|_\infty$ et on suppose que $\bx=(\frac{s}{g_2},\frac{2s}{g_2}, \ldots, \frac{(g_2-1)s}{g_2}, s)$ et $s\le R_\infty$.
Alors

\begin{eqnarray*}
\left|W_{2,\infty}(\bq)(1-|q|)^{\frac{g_2(g_2-1)}{2}}\prod_{k<l}(x_l-x_k)\right| &\ge& \left| W_{2,\infty}(\bq)\right|
\left(\frac{s(1-|q|)}{g_2}\right)^{\frac{g_2(g_2-1)}{2}}\\
&\ge& \left| W_{2,\infty}(\bq)\right| \left(\frac{0.995s}{g_2}\right)^{\frac{g_2(g_2-1)}{2}}.
\end{eqnarray*}

On choisissant un $\bq$ donn{\'e} par le lemme \ref{lemme:minorationduwronskien2} et en prenant
$s=\exp(-g^\cO)$ on prouve le

\begin{lemme}[Minoration du jacobien quadratique]\label{lemme:minorationdujacobien2}
Il existe une cons\-tan\-te effective $c_9$ telle que pour tout premier
$p$ tel que le genre  $g$ de $X_0(p)$ est au moins $2$, il existe  $\br= (r_1, \ldots, r_{g_2}) \in D(0,R_\infty)^{g_2}$
tel que $\log \left|\cJ_{2,\infty}(\br)\right| \ge -g^{c_9}.$ Un tel $\br$ se calcule en temps polynomial
en $p$ et la pr{\'e}cision absolue requise.
Les cinq normes sur $\cH^2(\Delta_2)$ d{\'e}finies pour $H=\sum_{\phi \in \cD_2} h_\phi \phi$ par $|H|=\max_\phi|h_\phi|$,
 $\cS_\infty(H)=\max_{q_\infty\in \bar
  D(0,R_\infty)}\left|\frac{q_\infty H}{(dq_\infty)^2}\right|$,
 $\cS_0(H)=\max_{q_0\in \bar D(0,R_0)}\left|\frac{q_0
  H}{(dq_0)^2}\right|$,  $\hat \cS_\infty (H)=\max_{q_\infty\in \bar
  D(0,\frac{1}{2})}\left|\frac{q_\infty H}{(dq_\infty )^2}\right|$,
et  $\hat \cS_0(H)=\max_{q_0\in \bar
  D(0,\frac{1}{2})}\left|\frac{q_0 H}{(dq_0)^2}\right|$, 
 sont dans des rapports $\exp(p^{\cO})$ ou, si l'on pr{\'e}f{\`e}re, 
il existe une constante positive  $c_4$  telle que pour tout $p$ et
pour tout $H\not = 0$  les 
diff{\'e}rences entre $\log |H|$, $\log \cS_0(H)$,  $\log \cS_\infty
(H)$, $\log \hat \cS_\infty  (H)$, et  $\log \hat \cS_0 (H)$  sont  born{\'e}es en valeur absolue
par $p^{c_4}$. \index{$\cS_\infty(H)$,  $\cS_0(H)$, 
$|H|$,  $\hat \cS_\infty(H)$,  et  $\hat \cS_0(H)$, les cinq  normes de la forme quadratique $H$}
\end{lemme}

\subsection{Stabilit{\'e}}\label{subsection:stabilite}

On suppose que $g\ge 2$.
On suppose choisie une origine $o$ et on note $S^gX$ la $g$-i{\`e}me puissance sym{\'e}trique 
de $X$ et $\mu^g : S^g X \rightarrow \CC^{\cD_1}/\cL$ l'application d'int{\'e}gration de Jacobi.
Soit $\epsilon = (\epsilon_k)_{1\le k \le g} \in \{ 0, \infty \}^g$ et $\br = (r_1, \ldots, r_g)\in 
D_\epsilon =  D(0,R_{\epsilon_1})\times \cdots  \times  D(0,R_{\epsilon_g})$ et notons $\rho \in \CC^g$ l'image de
$\br$ par $\mu_\epsilon$. On rappelle que l'instabilit{\'e} de $\br$ est l'oppos{\'e} du logarithme du module
du d{\'e}terminant jacobien $\cJ_{\cD_1,\epsilon}=\left| \frac{\omega}{dq_{\epsilon_k}}(r_{k})\right|_{\omega,k}$.
On suppose cette instabilit{\'e} finie et on la note $\lambda$. Donc  la restriction
de $\mu^g$ {\`a} un  voisinage  de $\br$ est injective. L'image d'un tel voisinage est
un voisinage de $\rho$ et  la corestriction de $\mu^g$ {\`a} un  voisinage assez petit
de $\rho$ est injective. On veut {\'e}tudier quantitativement cette
situation. On montre d'abord que les points voisins d'un point stable
sont assez stables. Ensuite on montre que l'image par l'application 
de Jacobi d'une boule centrée en un point stable, contient une boule
pas trop petite.

D'apr{\`e}s le lemme \ref{lemme:ordrejacobien}, le jacobien $\cJ_{\cD_1,\epsilon}$ 
est une s{\'e}rie en les $q_{\epsilon_k,k}$  d'ordre de grandeur $(g!, (3,\ldots,3))$.
Posons $\bx=(x_1, \ldots, x_g)$. La s{\'e}rie recentr{\'e}e
$\cJ_{\cD_1,\epsilon}(\br+\bx\star (\bun-\bB \br \bB ))$ est
d'ordre de grandeur $(\exp(g^\cO), (4, \ldots, 4))$. Pour
$|\bx|_\infty \le \frac{1}{2}$, le reste d'ordre $1$ de cette s{\'e}rie est major{\'e}
par $\exp(g^\cO)|\bx|_\infty $.
Donc si $|\bx|_\infty \le \exp(-g^\cO-\lambda)$ le point $\br+(\bun
-\bB \br \bB )\star \bx$ est stable d'instabilit{\'e} major{\'e}e
par $1+\lambda$. Si $\by= (\bun -\bB \br\bB )\star \bx$ alors $|\bx |_\infty \le p |\by|_\infty$ et $| \by|_\infty\le|\by|_2$
donc si $|\by|_2\le \exp(-g^\cO-\lambda)$ le point $\br+\by$  est stable d'instabilit{\'e} major{\'e}e
par $1+ \lambda$.

On compare maintenant l'application de Jacobi $\mu_\epsilon$ et sa diff{\'e}rentielle $D\mu_{\epsilon}$ en $\br$.
La matrice de $D\mu_{\epsilon}$ dans les bases
$(dq_{\epsilon_k})_k$ et $(z_\omega)_\omega$, la base duale de la base canonique
de $\CC^{\cD_1}$, n'est autre que $\left(  \frac{\omega}{dq_{\epsilon_k}}(r_{k})  \right)_{\omega, k}$. 
Elle a des coefficients  major{\'e}s par $3!2p^4$. Donc pour $\by\not = \bzero$ le quotient des
normes $L^\infty$ associ{\'e}es aux deux bases susmentionn{\'e}es $\frac{|D\mu_{\epsilon}\by|_\infty}{|\by|_\infty}$ est major{\'e}
par $g^\cO$.  Puisque le d{\'e}terminant de $D\mu_{\epsilon}$ 
est minor{\'e} en module par $\exp(-\lambda)$ on montre de même  que le quotient $\frac{|D\mu_{\epsilon}\by|_\infty}{|\by|_\infty}$
 est minor{\'e} par $\exp(-\lambda)g^{-\cO g}\ge \exp(-\lambda -\cO g^2)$ donc

\begin{equation}\label{equation:principale}
\exp(-\lambda -\cO g^2)|\by|_2\le |D\mu_{\epsilon}\by|_2 \le g^{\cO}|\by|_2.
\end{equation}

Soient $\omega$ une forme dans $\cD_1$ et $k$ un entier entre $1$ et $g$. Le d{\'e}veloppement de 
$\omega$ sur $D_{\epsilon_k}=D(0,R_{\epsilon_k})$ s'{\'e}crit $\omega
= \pm \frac{f(q_{\epsilon_k})}{q_{\epsilon_k}} d q_{\epsilon_k}$
pour une forme primitive, propre et normalis{\'e}e  $f$ de poids $2$. 
 La s{\'e}rie $f(q_{\epsilon_k})/q_{\epsilon_k}$ a pour ordre de grandeur
$(1,3)$. Sa recentr{\'e}e en $r_k$ a pour ordre de grandeur $(g^\cO,4)$.  Si $q_{\epsilon_k}=r_k+x_k(1-|r_k|)$,
avec $|x_k|\le \frac{1}{2}$, 
le reste d'ordre $1$ de la s{\'e}rie recentr{\'e}e est major{\'e} par $g^\cO|x_k|$.
Donc $\left|\int_{r_k}^{r_k+x_k(1-r_k)}\omega -x_k(1-r_k)\frac{f(r_k)}{r_k}\right|\le g^\cO (1-r_k)|x_k|^2$.
Ainsi le reste  d'ordre $2$

$$R_2{\mu_\epsilon}( (\bun -\br)\star \bx) = {\mu_\epsilon}(\br+(\bun
-\bB \br \bB )\star \bx)-{\mu_\epsilon}(\br)-D\mu_{\epsilon}((\bun -\bB \br \bB
)\star\bx )$$
\noindent  est major{\'e} en norme $L^\infty$
par $g^\cO|\bx|_\infty^2$ et en norme $L^2$ par $g^\cO|\bx|_2^2\le
g^{\cO}|\by|_2^2$ avec $\by=(\bun-\bB \br \bB )\star\bx$. Ainsi 

\begin{equation}\label{equation:erreur}
|R_2{\mu_\epsilon}( \by)|_2 = |{\mu_\epsilon}(\br+\by)-{\mu_\epsilon}(\br)-D\mu_{\epsilon}(\by )|_2\le g^{\cO}|\by|_2^2.
\end{equation}

Les {\'e}quations (\ref{equation:principale}) et (\ref{equation:erreur}) permettent de s'assurer
que la partie principale $D\mu_{\epsilon}\by$ est deux fois plus
grande que le reste $R_2\mu(\by)$ :  si 
$-\log |\by|_2 \ge \cO g^2+\lambda$ alors $|D\mu_{\epsilon}\by|_2\ge 2|R_2\mu( \by)|_2 $ donc
$\mu_\epsilon(\br+ \by)-\mu_\epsilon(\br)\ge \frac{1}{2}|D\mu_{\epsilon}\by|_2\ge \exp(-\lambda -\cO g^2)|\by|_2$.

Ainsi pour tout $\theta>  g^\cO$, la boule centr{\'e}e en $\rho$ et de rayon $\exp(-\cO g^2-2\lambda-\theta)$ 
pour la norme  $L^2$ est form{\'e}e
de points stables (ayant un ant{\'e}c{\'e}dent unique par l'application de Jacobi) et elle est contenue dans l'image
de la boule $L^2$ centr{\'e}e en $\br$ et de rayon $\exp(-\lambda-\theta)$.

\begin{lemme}[Stabilit{\'e} du probl{\`e}me inverse de Jacobi]\label{lemme:stabilitedejacobi}
Il existe une cons\-tante effective positive $c_{10}$   telle que l'{\'e}nonc{\'e} suivant soit vrai :

Soit $\epsilon = (\epsilon_k)_{1\le k \le g}$ dans $\{ 0, \infty \}^g$ et $\br = (r_1, \ldots, r_g)\in 
D_\epsilon =  D(0,R_{\epsilon_1})\times \cdots \times  D(0,R_{\epsilon_g})$ et notons $\rho$ l'image de
$\br$ par l'application d'int{\'e}gration de Jacobi  $\mu_\epsilon$. On suppose que l'instabilit{\'e} $\lambda$ de
$\br$ est finie. Alors pour tout $\theta \ge  g^{c_{10}}$ la boule $L^2$ centr{\'e}e en $\br$ et de rayon $\exp(-\lambda-\theta)$
est form{\'e}e de points stables et son image par l'application de Jacobi contient la boule centr{\'e}e en $\rho$ et
de rayon $\exp(-g^{c_{10}}-2\lambda-\theta)$. Donc  si $\br'\in S^gX$ s'envoie sur $\rho'$ par l'application de Jacobi
et si $|\rho-\rho'|_2\le \exp(-2g^{c_{10}}-2\lambda)$ alors $|\br'-\br|_2\le \exp(g^{c_{10}}+\lambda)|\rho-\rho'|_2$.
\end{lemme}

\subsection{Quelques sous ensembles discrets de la jacobienne}\label{subsection:reseaux}

On se place {\`a} nouveau dans la situation du paragraphe pr{\'e}c{\'e}dent. 
On suppose que $g\ge 2$.
On suppose choisie une origine $o$ et on note $S^gX$ la $g$-i{\`e}me puissance sym{\'e}trique 
de $X$ et $\mu^g : S^g X \rightarrow \CC^{\cD_1}/\cL$ l'application d'int{\'e}gration de Jacobi.
Soit $\epsilon = (\epsilon_k)_{1\le k  \le g} \in \{ 0, \infty \}^g$ et $\br = (r_1, \ldots, r_g)\in 
D_\epsilon =  D(0,R_{\epsilon_1})\times \cdots \times  D(0,R_{\epsilon_g})$ et notons $\rho$ l'image de
$\br$ par $\mu_\epsilon$. On note $\lambda$ l'instabilit{\'e} de $\br$. 

On identifie $\CC^g \supset D_\epsilon$ {\`a} l'espace tangent en $\br$ et on note $(\delta_k)_{1\le k\le g}$
la base duale de $(d q_{\epsilon_k})_{{1\le k\le g}}$. 

L'espace vectoriel r{\'e}el sous-jacent a pour base $(\delta_1,\delta_2, \ldots, \delta_g, i\delta_1, i\delta_2,\ldots,i\delta_g)$.
Pour $k$ entre $1$ et $g$ on pose $i\delta_k=\delta_{k+g}$.

On choisit un r{\'e}el $\chi > \log p$,
on pose $\Upsilon = \exp(-\chi)$  et on consid{\`e}re
les $2g$ petits accroissements $\beta_1$, \ldots, $\beta_{2g}$ de $\CC^g/\cL$  d{\'e}finis par
$\beta_k=\mu(r_k+\Upsilon)-\mu(r_k)=(\int_{r_k}^{r_k+\Upsilon}\omega)_{\omega
  \in \cD_1}= \mu_\epsilon (\br +\Upsilon \delta_k)
-\mu_\epsilon (\br)$ pour $1\le k \le g$ et
$\beta_{g+k}=\mu(r_k+i\Upsilon)-\mu(r_k)=(\int_{r_k}^{r_k+i\Upsilon}\omega)_{\omega
  \in \cD_1}= \mu_\epsilon (\br +\Upsilon i\delta_k)
-\mu_\epsilon (\br)$ pour $1\le k \le g$.

Si $M$ est un entier positif, on note $\cA(\br,\chi,M)$\index{$\cA(\br,\chi,M)$, sous ensemble discret 
du tore complexe} le sous ensemble de $\CC^g/\cL$ form{\'e}
des combinaisons des $\beta_k$ {\`a} coefficients entiers dans $[-M,M]$.
On se demande dans quelle mesure les {\'e}l{\'e}ments du tore $\CC^g/\cL$ peuvent {\^e}tre approch{\'e}s par
l'ensemble $\cA(\br,\chi,M)$.

On doit premi{\`e}rement minorer le d{\'e}terminant $\frac{\beta_1\wedge \cdots \wedge \beta_{2g}}{e_1\wedge\cdots\wedge e_{2g}}$
o{\`u} $(e_k)_k$ est la base canonique de l'espace r{\'e}el $\RR^{2g}$ sous-jacent {\`a} $\CC^g=\CC^{\cD_1}$.  En particulier
$(e_{k})_{1\le k \le g}$  est la base canonique de $\CC^g$ et $e_{k+g}=ie_k$.
On note $D\mu_{\epsilon}$ la diff{\'e}rentielle de $\mu_\epsilon$ en $\br$. Pour $1\le k\le 2g$ on note
    $\gamma_k=D\mu_{\epsilon}(\delta_k)$.
Le d{\'e}terminant $\frac{\gamma_1 \wedge \cdots \wedge \gamma_{2g}}{e_1\wedge \cdots \wedge e_{2g}}$ est le carr{\'e} du module
du d{\'e}terminant jacobien.

Or $\Upsilon \gamma_k$ est une bonne approximation de $\beta_k$.
Plus précisément, on a vu au   paragraphe  pr{\'e}c{\'e}dent que si $-\log \Upsilon =
\chi \ge  \cO g^2
+\lambda$ alors $\Upsilon |\gamma_k|_2\ge 2|\beta_k-\Upsilon \gamma_k|_2$
donc $\frac{\Upsilon}{2} |\gamma_k|_2\le  |\beta_k|_2\le \frac{3\Upsilon}{2}|\gamma_k|_2\le g^{\cO}\Upsilon$
d'apr{\`e}s l'in{\'e}galit{\'e} (\ref{equation:principale}).

L'in{\'e}galit{\'e} (\ref{equation:erreur})  implique quant {\`a} elle 
que $|\beta_k-\Upsilon \gamma_k|_2 \le g^{\cO}\Upsilon ^2$.

Comme le d{\'e}terminant est multilin{\'e}aire, 
on majore le module de la diff{\'e}rence $\frac{\beta_1\wedge \cdots \wedge \beta_{2g}}{e_1\wedge\cdots\wedge e_{2g}}
-\Upsilon^{2g}\frac{\gamma_1 \wedge \cdots \wedge \gamma_{2g}}{e_1\wedge \cdots \wedge e_{2g}}$ par $2^{2g}\left(\max_{1\le k\le 2g}|\beta_k|_2\right)^{2g-1}
\max_{1\le k\le 2g}|\beta_k-\Upsilon\gamma_k|_2$.

Ainsi la diff{\'e}rence $\frac{\beta_1\wedge \cdots \wedge \beta_{2g}}{e_1\wedge\cdots\wedge e_{2g}}
-\Upsilon^{2g}\frac{\gamma_1 \wedge \cdots \wedge \gamma_{2g}}{e_1\wedge \cdots \wedge e_{2g}}$ est major{\'e}e 
par $2^{2g}g^{\cO(2g-1)}g^{\cO}\Upsilon^{2g+1}$.
Cette diff{\'e}rence est plus petite que la moiti{\'e} de $\Upsilon^{2g}|\frac{\gamma_1 \wedge \cdots \wedge \gamma_{2g}}{e_1\wedge \cdots \wedge e_{2g}}|=
\exp(-2\lambda)\Upsilon^{2g}$ pourvu que $\chi \ge  \cO g^2 +2\lambda$.   Alors 

\begin{equation*}
\frac{\beta_1\wedge \cdots \wedge \beta_{2g}}{e_1\wedge\cdots\wedge e_{2g}} \ge \frac{1}{2}\exp(-2\lambda-2g\chi).
\end{equation*}

On a donc une minoration du d{\'e}terminant de la matrice de passage de la base $(e_k)_{1\le k \le 2g}$ {\`a} la
base $(\beta_k)_{1\le k \le 2g}$. 
Les colonnes de cette matrice sont  les $\beta_k$ et elles sont  major{\'e}s en   norme $L^\infty$ par $g^{\cO} \Upsilon$.
Donc les colonnes de la matrice inverse sont major{\'e}es en norme $L^\infty$ par $2(2g-1)!\exp(2\lambda+2g\chi)g^{\cO(2g-1)}\Upsilon^{2g-1}$.

D'apr{\`e}s le lemme \ref{lemme:lereseaudesperiodes} il existe une constante positive 
$c_2$ et une base du r{\'e}seau $\cL$ des p{\'e}riodes de coordonn{\'e}es major{\'e}es en valeur absolue par $\exp(p^{c_2})$ dans
la base canonique de $\RR^{2g}$. Donc les coordonn{\'e}es de ces p{\'e}riodes dans la base $(\beta_k)_{1\le k\le 2g}$ sont born{\'e}es 
par $2(2g-1)!\exp(p^{c_2}+2\lambda+2g\chi)g^{\cO(2g-1)}\Upsilon^{2g-1} \le \exp(2\lambda+\chi+g^{\cO} )$.

Les coordonnées de tous les points du parall{\'e}logramme  fondamental associ{\'e} {\`a}
cette base sont également bornés de la sorte.
Enfin, en rempla{\c c}ant de telles coordonn{\'e}es par leurs valeurs enti{\`e}res on commet une erreur 
born{\'e}e en norme $L^2$ dans la base canonique par $g^{\cO}\Upsilon$.  

\begin{lemme}[Bonne r{\'e}partition]\label{lemme:bonnerepartition}
Il existe une constante $c_5$ effective positive telle que :  si $\epsilon \in \{0, \infty \}^g$ et 
$\br \in D_\epsilon$ a une  instabilit{\'e} finie $\lambda$ et si 
 $\chi$ est un r{\'e}el  plus grand que $g^{c_5}+2\lambda$ et si   $M$ est le plus
petit  entier plus grand que $\exp( g^{c_5} + 2\lambda 
+\chi)$, alors tout point $\rho$ de $\CC^g=\RR^{2g}$ est \`a distance $\le \exp(-\chi+g^{c_5})$ d'un point 
de $\cA(\br,\chi, M)$ et les coefficients entiers dans $[-M,M]$ de ce 
 point peuvent {\^e}tre calcul{\'e}s en temps polynomial en $g$, $\log (\max(1,|\rho|_\infty))$ et $\chi$.
\end{lemme}

\section{Complexit{\'e} des op{\'e}rations dans la jacobienne}\label{section:arithmetique}

Dans le paragraphe \ref{subsection:dualisation} on {\'e}tudie la complexit{\'e} de l'algorithme qui, {\'e}tant donn{\'e} une famille
de points, trouve une forme s'annulant en ces points, puis recherche les autres z{\'e}ros de cette forme.
On en d{\'e}duit dans le paragraphe \ref{subsection:additionsoustraction} que les op{\'e}rations {\'e}l{\'e}mentaires
dans $J_0(p)$ sont stables et se font en temps déterministe polynomial en $p$ et la pr{\'e}cision requise.

Au paragraphe \ref{subsection:inversejacobi}, un semblable r{\'e}sultat est {\'e}tabli pour le probl{\`e}me inverse de Jacobi.

\subsection{Une op{\'e}ration {\'e}l{\'e}mentaire  dans la jacobienne : la dualisation }\label{subsection:dualisation}

Un point $P$  de $X$ est sp{\'e}cifi{\'e} par une valeur de $\tau \in \cH^*$ ou de 
$q_0\in  D(0,R_0)$ ou de $q_\infty  \in D(0,R_\infty)$.  Conna{\^\i}tre un tel point, c'est disposer
d'un algorithme qui calcule par exemple  le $q_\infty$  correspondant, en temps polynomial
en $g$ et la pr{\'e}cision absolue requise.

On suppose $g\ge 2$ et on se donne $3g-4$ points distincts   $R_1$, $R_2$, \ldots, $R_{3g-4}$ sur $X$. On suppose
que ces points ne sont ni des pointes, ni des points elliptiques. 
On pose $e=3g-2+\nu_2+\nu_3=g_2-1$ et on note $P_1$, \ldots, $P_{e}$ la famille de points form{\'e}e des $3g-4$ points
que l'on vient de se donner et des $2+\nu_2+\nu_3$ points elliptiques ou paraboliques.

On suppose  que les $r$ premiers
points $P_1$, \ldots, $P_r$ sont dans $D_\infty$ et sont donc donn{\'e}s par des $q_{\infty,k}\in D(0, R_\infty)$ pour
$1\le k\le r$. Les
autres points sont suppos{\'e}s appartenir {\`a} $D_0$ et sont donc donn{\'e}s par des $q_{0,k}\in D(0, R_0)$
pour $r+1\le k\le e$.  On suppose que $P_1$ est la pointe $\infty$ et $P_e$ la pointe $0$.

On se donne un $\eta >0$ et pour tout $k$ entre $1$
et $r$ on consid{\`e}re le disque $D_{\infty, k}=\{q_\infty, |q_\infty - q_{\infty,k}|\le \eta \}$ de rayon $\eta$
centr{\'e} en $q_{\infty,k}$. De m{\^e}me, pour tout $k$ entre $r+1$ et $e$ on note 
$D_{0, k}=\{q_0, |q_0 - q_{0,k}|\le \eta \}$ le disque de rayon $\eta$
centr{\'e} en $q_{0,k}$. On suppose que les images de  tous ces disques dans $X(\CC)$
sont deux {\`a} deux disjointes.

\begin{definition}[Résolution d'une famille de points]
Dans cette situation, on dit  que $((P_k)_{1\le k \le 3g-2+\nu_2+\nu_3}, r, \eta)$ 
est une {\it r{\'e}solution} pour la famille de points 
$(R_k)_{1\le k \le 3g-4}$.
\end{definition}\index{r{\'e}solution d'une famille de points}

On note $P=P_1+\cdots+P_{e}$ le diviseur somme des $P_k$ et  on observe que $\cH^2(-P+\Delta_2) \subset
\cH^2(\Delta_2)$ est un  sous-espace vectoriel non nul de $\cH^2$.  

On repr{\'e}sente une forme $F$ dans  $\cH^2(\Delta_2)$ par ses coordonn{\'e}es dans la base $\cD_2$ soit
$F=\sum_{\phi \in \cD_2} f_\phi \phi$. On cherche les {\'e}quations lin{\'e}aires qui d{\'e}finissent 
le sous espace  $\cH^2(-P+\Delta_2)$ dans cette base. Il y a une {\'e}quation pour chaque $P_k$.

L'{\'e}quation correspondant {\`a} la pointe {\`a} l'infini $P_1$,  est 

$$\sum_{\phi \in \cD_2}\frac{q_\infty\phi }{(dq_\infty)^2}(P_1)f_\phi = \sum_\phi a_{\phi, 1}f_\phi = 0$$
\noindent o{\`u} $\frac{q^2_\infty\phi }{(dq_\infty)^2}=\sum_{k\ge 1}a_{\phi, k}q_\infty^k$ est le d{\'e}veloppement
en l'infini de la forme primitive, propre et normalis\'ee  de poids $4$ associ{\'e}e {\`a} $\phi$.

L'{\'e}quation correspondant {\`a} la pointe en z{\'e}ro  $P_e$,  est 

$$\sum_{\phi \in \cD_2}\frac{q_0\phi }{(dq_0)^2}(P_e)f_\phi = \sum_\phi b_{\phi, 1}f_\phi = 0$$
\noindent o{\`u} $\frac{q^2_0\phi }{(dq_0)^2}=\sum_{k\ge 1}b_{\phi, k}q_0^k$ est le d{\'e}veloppement
en z{\'e}ro de la forme primitive, propre et normalis\'ee  de poids $4$ associ{\'e}e {\`a} $\phi$. Noter que $b_{\phi,1}=\pm a_{\phi, 1}$ car les
formes primitives, propres et normalis{\'e}es  sont en particulier valeurs propres de l'involution d'Atkin-Lehner.

L'{\'e}quation correpondant {\`a} $P_k$ pour $2\le k\le r$ est 
$$\sum_{\phi \in \cD_2}\frac{q_\infty\phi }{(dq_\infty)^2}(P_k)f_\phi =  0$$
\noindent o{\`u} $\frac{q_\infty\phi }{(dq_\infty)^2}(P_k)=\sum_{l\ge 1}a_{\phi, l}q_{\infty, k}^{l-1}$ 
se calcule en temps polynomial en $p$ et la pr{\'e}cision requise. 

L'{\'e}quation correpondant {\`a} $P_k$ pour $r+1\le k\le e-1$ est 
$$\sum_{\phi \in \cD_2}\frac{q_0\phi }{(dq_0)^2}(P_k)f_\phi =  0$$
\noindent o{\`u} $\frac{q_0\phi }{(dq_0)^2}(P_k)=\sum_{l\ge 1}b_{\phi, l}q_{0, k}^{l-1}$ 
se calcule en temps polynomial en $p$ et la pr{\'e}cision requise.

On rassemble toutes ces {\'e}quation dans une matrice $\cM_P$ dont le noyau
d{\'e}crit $\cH^2(-P+\Delta_2)$ dans la base $\cD_2$. Les coefficients de
cette matrice sont $\le p^{\cO}$. On calcule une valeur approch{\'e}e
$\cM_P'$ de $\cM_P$ {\`a} coefficients dans $\QQ (i)$. On note $m\ge 2$ la pr{\'e}cision
de cette approximation. Donc $\cM_P'-\cM_P$ a des coefficients 
inf{\'e}rieurs {\`a} $\exp(-m)$ en module. Il existe une telle $\cM_P'$ de taille
polynomiale en $\log p$ et $m$ et on la calcule en temps polynomial
en $p$ et $m$. L'algorithme d'Hermite \cite{knuth}, ou mieux l'algorithme LLL \cite[Theorem 2.6.2]{cohen},  produisent 
une base {\`a} coefficients entiers du noyau de $\cM_P'$ de taille polynomiale en $p$ et $m$ en temps
polynomial en $p$ et $m$.
On choisit  un vecteur de cette base. Si $M=a+ib$ est le plus grand des coefficients de ce vecteur,
on normalise en divisant le vecteur par le maximum de $|a|$ et $|b|$. On obtient ainsi 
un  $f'=(f'_\phi)_{\phi \in \cD_2}$ dans le noyau de $\cM_P'$ {\`a} coefficients dans $\QQ(i)$ et 
   $F'=\sum_\phi f'_\phi \phi$ la forme associ{\'e}e telle que $|F|=\max_\phi |f'_\phi|$ soit compris entre $1$
et $\sqrt 2$. 
D'apr{\`e}s le lemme  \ref{lemme:minorationdujacobien2}, les normes
$\hat \cS_0(F')$ et $\hat \cS_\infty(F')$ sont 
major{\'e}es par $\sqrt 2 \exp(p^{c_4})$ et minor{\'e}es par $\exp(-p^{c_4})$
pour une certaine constante positive $c_4$.

Le produit $\cM_P f' = (\cM_P - \cM'_P)f'$ est un vecteur de coefficients major{\'e}s en module par  $g_2\exp(-m)\sqrt 2$. Ces coefficients sont, dans l'ordre,
la valeur de $\frac{q_\infty F'}{(dq_\infty)^2}$ en la pointe
infini, les valeurs de $\frac{q_\infty F'}{(dq_\infty)^2}$ en
les $P_k$ pour $2\le k \le r$, les valeurs de 
$\frac{q_0 F'}{(dq_0)^2}$ en
les $P_k$ pour $r+1\le k \le e-1$, la valeur 
de 
$\frac{q_0 F'}{(dq_0)^2}$ en
la pointe nulle $P_e$.

On peut maintenant appliquer le lemme \ref{lemme:stabunzero} de stabilit{\'e} des z{\'e}ros. 

On s'int{\'e}resse d'abord aux $P_k$ pour $k\le r$. Ils sont contenus dans $D(0,R_\infty)$. Puisque $|F'|\le \sqrt 2$,
la s{\'e}rie enti{\`e}re $\frac{q_\infty F'}{(dq_\infty)^2}$ est d'ordre de grandeur $(g_2\sqrt 2,4)$. Avec les notations
du lemme \ref{lemme:stabunzero} on a donc $A=g_2\sqrt 2$, $n=4$,  $1-|c|\ge 0.995$ et $-\log \epsilon = m/2$ pourvu que
$m\ge \cO \log (g)$. On demande que $\frac{(0.995)^{14}\sqrt\frac{m}{2}}{c_{14}}$ soit plus grand
que $c_{14}( 4^2(1-\log {0.995})+\log g_2\sqrt 2)$ et que $p^{c_4}$ ce qui est assur{\'e} si
$m\ge p^{\cO}$. Alors chaque $P_k$ pour $k\le r$ est {\`a} distance $\le \exp(-\sqrt{m/2})$ d'un z{\'e}ro de $F'$.

On s'int{\'e}resse maintenant aux $P_k$ pour $k>r$. Ils sont contenus dans $D(0,R_0)$ avec
$R_0=1-\frac{1}{p}$. Puisque $|F'|\le \sqrt 2$,
la s{\'e}rie enti{\`e}re $\frac{q_0 F'}{(dq_0)^2}$ est d'ordre de grandeur $(g_2\sqrt 2,4)$. Avec les notations
du lemme \ref{lemme:stabunzero} on a donc $A=g_2\sqrt 2$, $n=4$, $1-|c|\ge \frac{1}{p}$
    et $-\log \epsilon = m/2$ pourvu que
$m\ge \cO \log (g)$. 
On demande que $\frac{\sqrt\frac{m}{2}}{p^{14}c_{14}}$ soit plus grand
que $c_{14}( 4^2(1+\log {p})+\log g_2\sqrt 2)$ et que $p^{c_4}$ ce qui est assur{\'e} si
$m\ge p^{\cO}$. 
 Alors chaque $P_k$ pour $k> r$ est {\`a} distance $\le \exp(-\sqrt{m/2})$ d'un z{\'e}ro de $F'$.

On s'assure enfin que $\exp(-\sqrt{m/2})$ est plus petit que la r{\'e}solution $\eta$ de sorte
que les $e$ z{\'e}ros ainsi trouv{\'e}s sont distincts. On cherche alors les $g$  autres z{\'e}ros de $F'$ avec l'algorithme
de quadrichotomie de Weyl (voir le lemme \ref{lemme:stabglobalezeros}). Puisque l'application d'int{\'e}gration de Jacobi
est Lipschitzienne avec coefficients de dilatation polynomiaux en $p$  on obtient le

\begin{lemme}[Dualisation]
Il existe un algorithme d{\'e}terministe qui, {\'e}tant don\-n{\'e}s un nombre premier $p$ et 
$3g-4$ points distincts $R_1$, \ldots, $R_{3g-4}$ sur $X_0(p)$
et une  r{\'e}solution $((P_m)_{1\le m < g_2},r,\eta )$, calcule $g$ points $Q_1$, \ldots, $Q_g$ dont la somme
est lin{\'e}airement {\'e}quivalente {\`a} deux fois le diviseur canonique moins la somme des points initiaux :

$$Q_1+\cdots+Q_{g}\sim 2\cK - (R_1+\cdots+R_{3g-4}).$$

Cet algorithme est polynomial en $p$, le logarithme $-\log \eta$ de la r{\'e}solution et la pr{\'e}cision requise.
\end{lemme}

{\bf Remarque importante}: Notons   que dans cet {\'e}nonc{\'e}, l'erreur qui est major{\'e}e se
mesure dans l'image $\CC^{\cD_1}$ de l'application
de Jacobi.

On peut s'affranchir de la condition sur la r{\'e}solution $\eta$. En effet, si les points $R_k$ ne sont pas
tous disctincts, ou plus g{\'e}n{\'e}ralement  si leur r{\'e}solution est jug{\'e}e trop faible, 
c'est {\`a} dire si $\eta$ est trop petit, on peut perturber l{\'e}g{\`e}rement ces points
pour obtenir la r{\'e}soluton souhait{\'e}e. Comme l'application d'int{\'e}gration de Jacobi est Lipschitzienne
de dilatation polynomiale en $p$, on peut
adapter la perturbation {\`a} la pr{\'e}cision finale souhait{\'e}e. On obtient ainsi le 

\begin{lemme}[Dualisation]\label{lemme:dualisation}
Il existe un algorithme d{\'e}terministe qui, {\'e}tant don\-n{\'e}s un nombre premier $p$ et 
$3g-4$ points  $R_1$, \ldots, $R_{3g-4}$ sur $X_0(p)$, calcule $g$ points $Q_1$, \ldots, $Q_g$ dont la somme
est lin{\'e}airement {\'e}quivalente {\`a} deux fois le diviseur canonique moins la somme des points initiaux :

$$Q_1+\cdots+Q_{g}\sim 2\cK - (R_1+\cdots+R_{3g-4}).$$

Cet algorithme est polynomial en $p$ et la pr{\'e}cision requise.
\end{lemme}

\subsection{Addition et soustraction dans la jacobienne}\label{subsection:additionsoustraction}

On suppose que le genre de $X_0(p)$ est au moins $4$.
On choisit  un diviseur effectif origine de degr{\'e} $g$ not{\'e} $O=o_1+\cdots+o_g$.
On choisit aussi un diviseur effectif auxiliaire de degr{\'e} $g-4$ not{\'e} $N$.
La dimension de $\cH^2(\Delta_2)$  est $g_2=3g-1+\nu_2+\nu_3$. 
Le plus naturel est de choisir un point $o$ origine de l'application d'int{\'e}gration de Jacobi
et de poser $O=go$ et $N=(g-4)o$. 
Un {\'e}l{\'e}ment de $\Pic^0(X)$ est donn{\'e} comme classe d'un diviseur $Q-O$ o{\`u} $Q$
est effectif de degr{\'e} $g$.  Soit  $R$ un autre diviseur effectif de degr{\'e} $g$. Pour ajouter la classe de $Q-O$  et celle  
de $R-O$ on applique le lemme \ref{lemme:dualisation}
au diviseur $Q+R+N$ qui est effectif de degr{\'e} $3g-4$.
On obtient un diviseur $T$ effectif de degr{\'e} $g$ tel que $T\sim 2\cK-Q-R-N$.
On applique {\`a} nouveau le lemme \ref{lemme:dualisation} au diviseur $T+O+N$ et on obtient
un diviseur $U$ effectif  de degr{\'e} $g$ tel que $U+O\sim Q+R$. Donc $U-O$ est bien
la somme de $Q-O$ et $R-O$.

Pour calculer l'oppos{\'e} de $Q-O$ on applique le lemme \ref{lemme:dualisation} au diviseur
$2O+N$ ce qui produit un diviseur $\Xi$ effectif de degr{\'e} $g$  {\'e}quivalent {\`a} $2\cK-N-2O$.
On applique {\`a} nouveau le lemme \ref{lemme:dualisation} au diviseur $\Xi+Q+N$ et on trouve un diviseur
$V$ effectif de degr{\'e} $g$ tel que $V-O\sim -(Q-O)$.

\begin{theoreme}[Arithm{\'e}tique dans la jacobienne]\label{theoreme:arithmetiquedebase}
Les op{\'e}ration d'addition et de soustraction  dans la jacobienne de $X_0(p)$ se font en temps
d{\'e}terministe polynomial en $p$ et la pr{\'e}cision requise.
\end{theoreme}

\subsection{Le probl{\`e}me inverse de  Jacobi}\label{subsection:inversejacobi}

\'Etant donn{\'e}  $\rho=(\rho_\omega)_\omega \in \CC^{\cD_1}$ on cherche un diviseur effectif
de degr{\'e} $g$ sur $X$ qui s'envoie sur $\rho$ modulo le r{\'e}seau $\cL$ des p{\'e}riodes par l'application de Jacobi.

Le lemme \ref{lemme:bonnerepartition} de bonne r{\'e}partition, que l'on applique au $\br$ fourni
par le lemme \ref{lemme:minorationdujacobien} de minoration du jacobien, et {\`a} un r{\'e}el $\chi$
assez grand, fournit un ensemble $\cA(\br, \chi,M)=\{t_1 \beta_1+\cdots+t_{2g}\beta_{2g} \mbox{ avec }  -M\le t_k \le M   \}$
o{\`u} les $\beta_k$ sont les images par l'application de Jacobi de diviseurs $\Omega_k$ connus, diff{\'e}rences de deux points
proches. Les coordonn{\'e}es r{\'e}elles de $\rho$ dans la base form{\'e}e des $\beta_k$ sont calcul{\'e}es en inversant
une matrice dont le d{\'e}terminant n'est pas trop petit, puis tronqu{\'e}es {\`a} l'entier inf{\'e}rieur ou {\'e}gal
le plus proche.  Les entiers $t_k$ ainsi  obtenus, on calcule dans $J_0(p)$ la combinaison
$\sum_{1\le i \le 2g}t_k\Omega_k$ avec la m{\'e}thode du  th{\'e}or{\`e}me \ref{theoreme:arithmetiquedebase} d'arithm{\'e}tique
dans la jacobienne et en utilisant l'exponentiation rapide. Le nombre d'op{\'e}rations {\'e}l{\'e}mentaires dans la
jacobienne 
est donc polynomial en $\log M$, donc aussi la perte de pr{\'e}cision. 

\begin{theoreme}[Problème inverse de Jacobi]\label{theoreme:inverse}
\'Etant donn{\'e}  un entier  premier $p$,  on note $g$ le genre de $X_0(p)$ et
$\cD_1$
 l'ensemble des formes diff{\'e}rentielles primitives, propres et normalis{\'e}es  et $\cL$ le r{\'e}seau des p{\'e}riodes
de $X_0(p)$.

Il existe un algorithme d{\'e}terministe qui pour  $\rho \in 
\CC^{\cD_1}$  et pour $O=o_1+\cdots+o_g$ une origine de l'application
d'int{\'e}gration de Jacobi $\mu ^g : S^g X_0(p) \rightarrow \CC^{\cD_1}/\cL$ 
et pour  $k$ un entier positif, calcule $g$ points $P_1$, \ldots, $P_g$ de $X_0(p)$ tels
que $\left| \mu^g(P_1, \ldots, P_g) -\rho \right|_\infty \le \exp(-k)$,  en temps polynomial
en $p$,    la pr{\'e}cision $k$, et  la taille $\log (\max(1,|\rho|_\infty))$ de $\rho$.
\end{theoreme}

On note que la d{\'e}pendance en  $\log (\max(1,|\rho|_\infty))$ est sans importance car le r{\'e}seau
des p{\'e}riodes admet un parall{\'e}logramme fondamental de rayon  $\exp(p^{\cO})$ d'apr{\`e}s le lemme
\ref{lemme:lereseaudesperiodes}.

Notons  encore que dans cet {\'e}nonc{\'e}, l'erreur qui est major{\'e}e se mesure dans l'image de l'application
de Jacobi. Si l'on veut majorer l'erreur commise sur les points $P_k$ on note
$Q_1$, \ldots, $Q_g$ des points tels que $\mu^g(Q_1,\ldots,Q_g)=\rho$ et soit $\lambda$
leur instabilit{\'e} que l'on suppose finie, sans quoi la question serait d{\'e}pourvue de sens.

Le lemme \ref{lemme:stabilitedejacobi} de stabilit{\'e} du probl{\`e}me inverse de Jacobi, montre
que la perte de pr{\'e}cision dans le probl{\`e}me inverse de Jacobi est polynomiale en $p$ et $\lambda$.
Il s'agit donc de controler $\lambda$.

Soit $\epsilon \in \{0, \infty \}^g$ tel que $Q_k\in D_{\epsilon_k}$ pour tout $1\le k \le g$ et notons
$q_{\epsilon_k, k}$ la valeur de $q_{\epsilon_k}$ en $Q_k$. L'instabilit{\'e}  $\lambda$  est l'inverse du logarithme du
module du d{\'e}terminant jacobien $\left| \frac{\omega}{dq_{\epsilon_k}}(q_{\epsilon_k,k})\right|_{\omega,k}$
qui n'est pas une fontion alg{\'e}brique.
Il est donc naturel de r{\'e}{\'e}crire ce jacobien comme le produit d'une quantit{\'e} alg{\'e}brique et 
de facteurs plus simples.

Pour tout $k$ entre $1$ et $g$, on suppose que  la valeur $j(Q_k)$ de l'invariant modulaire $j$
en $Q_k$ v{\'e}rifie $j(Q_k)\not\in \{0, 1728, \infty \}$.

On note que  $j(Q_k)$   ne peut {\^e}tre proche de 
$0$ et $1728$ en m{\^e}me temps. Si $j(Q_k)$ n'est pas proche de $0$ on pose $j_k = j-1728$ et sinon on pose
$j_k = j$.  Le jacobien se r{\'e}{\'e}crit alors comme produit du d{\'e}terminant 
$\cJ_{alg}=
\left| \frac{j_k\omega}{dj}(q_{\epsilon_k,k})\right|_{\omega,k}$ et  des $\frac{dj}{j_kdq_{\epsilon_k}}(q_{\epsilon_k,k})$
pour $1\le k \le g$.

Supposons que les $Q_k$ sont des points alg{\'e}briques sur $\bar \QQ$ et  que l'instabilit{\'e} est finie. Alors $\cJ_{alg}$ est fini et  non nul et on peut le minorer
en fonction de la hauteur des $j(Q_k)$. Les $\frac{dj}{j_kdq_{\epsilon_k}}(q_{\epsilon_k,k})$ se minorent
en utilisant le lemme \ref{lemme:stabglobalezeros} de stabilit{\'e} des z{\'e}ros d'une s{\'e}rie enti{\`e}re\footnote{ou
un argument de compacit{\'e}, plus g{\'e}n{\'e}ral mais non effectif.} : si 
$\frac{dj}{dq_{\epsilon_k}}(q_{\epsilon_k,k})$ est petit alors $q_{\epsilon_k,k}$
est proche d'un z{\'e}ro de $\frac{dj}{ dq_{\epsilon_k} }$. Mais  le
 $j_k$ choisi est alors plus petit encore.

\appendix

\section{Appendice  sur les s{\'e}ries enti{\`e}res}\label{section:analyse}

On a rassembl{\'e} dans cette section  les  d{\'e}finition et r{\'e}sultats relatifs {\`a} la localisation
et {\`a} la stabilit{\'e} des z{\'e}ros des s{\'e}ries enti{\`e}res qui sont n{\'e}cessaires {\`a} notre travail. 

Le  paragraphe \ref{subsection:ordredegrandeur} introduit quelques d{\'e}finitions et notations ainsi que des r{\'e}sultats {\'e}l{\'e}mentaires.
Le paragraphe \ref{subsection:prol} {\'e}nonce et d{\'e}montre une forme simple et quantifi{\'e}e du th{\'e}or{\`e}me de prolongement
analytique.  On introduit le polygone de Newton d'une s\'erie enti\`ere
dans
le paragraphe \ref{subsection:newton} et on montre comment il permet
de localiser les z\'eros de cette s\'erie. Le paragraphe \ref{subsection:stabilitezeros} prouve un r{\'e}sultat de stabilit{\'e} des z{\'e}ros d'une s{\'e}rie enti{\`e}re
et en d{\'e}duit une majoration de la complexit{\'e} de la localisation des z{\'e}ros.
 
\subsection{Ordre de grandeur}\label{subsection:ordredegrandeur}

Soit $g\ge 1$ un entier.
Si $\bx=(x_1,\ldots,x_g)\in \CC^g$, on note $|\bx |_\infty=\max_k|x_k|$  la norme
$L^\infty$ de $\bx$. On note $|\bx|_1=\sum_k|x_k|$ la norme $L^1$ de
$\bx$ et $|\bx|_2=\sqrt{\sum_k|x_k|^2}$ sa norme $L^2$.
On note $\bB \bx \bB$ le vecteur $(|x_1|, \ldots, |x_g|)$. Si 
$\by=(y_1,\ldots,y_g)$ on note $\bx \star \by$\index{$\star$} le vecteur
$(x_1y_1,\ldots, x_gy_g)$.
On note $\bzero=(0,\ldots,0)\in \CC^g$ et $\bun=(1,\ldots,1)
\in \CC^g$. Si $\bx=(x_1,\ldots,x_g)\in \RR^g$, on dit que $\bx \ge \bzero$ si et seulement
si $x_k \ge 0$ pour tout $k$. On dit que $\bx > \bzero$ si et seulement
si $x_k > 0$ pour tout $k$.
Notons $P(\bx, \br)=\prod_{k=1}^g D(x_k, r_k) \subset \CC^g$\index{$P(\bx,\br)$, le polydisque
de centre $\bx$ et de polyrayon $\br$} le polydisque de polycentre $\bx$
 et de polyrayon $\br$.

Une s{\'e}rie enti{\`e}re $f$ est donn{\'e}e par  $f=\sum_{\bk}f_\bk\bx^\bk$ o{\`u}
$\bk$ parcourt $\NN^g$.

\begin{definition}[Ordre de grandeur]\index{ordre de grandeur d'une s{\'e}rie}
Soit $A\ge 1$ un r{\'e}el et $\bn =(n_1,\ldots,n_g)\in \NN^g$ tel que
$\bn \ge \bun$. On dit que
$f$ est d'ordre de grandeur $(A, \bn)$ si  pour tout $\bk \ge \bzero$
on a 

$$|f_\bk|\le A(\bk+\bun)^\bn=A\prod _{1\le m \le g}(k_m+1)^{n_m}.$$

Si $f$ est d'ordre de grandeur $(A, \ba)$ et $h$ d'ordre de grandeur
$(B, \bb)$ alors le produit $p=fh$ est d'ordre de grandeur $(AB,\ba+\bb+\bun)$.
\end{definition}

En effet $p=\sum_\bm p_\bm\bx^\bm$ avec $p_\bm=\sum_{\bk+\bl=\bm}f_\bk h_\bl$. Il
y  a $\prod_{1\le n\le g}(m_n+1)=(\bm+\bun) ^\bun$ termes dans cette derni{\`e}re somme, et chacun
est major{\'e} en module par  $|f_\bk h_\bl|\le A(\bk+\bun)^\ba B(\bl+\bun)^\bb\le AB(\bk+\bl+\bun)^{\ba+\bb}$.\hfill $\Box$

\begin{lemme}[D{\'e}riv{\'e}e]\label{lemme:derivprod}
Si $f$ est une s{\'e}rie enti{\`e}re  d'une variable, d'ordre de grandeur $(A,m)$, alors sa d{\'e}riv{\'e}e $d$-i{\`e}me est
d'ordre de grandeur $(A2^{dm+\frac{d(d-1)}{2}},m+d)$. 
\end{lemme}

Soit  $g\ge 1$ un entier et $f=\sum_\bk f_\bk\bz^\bk$ une s{\'e}rie enti{\`e}re de $g$ variables d'ordre de grandeur
 $(A, \bn)$. On en d{\'e}duit pour
tout $\bz\in P(\bzero,1)$ la majoration

\begin{eqnarray*}
|f(\bz)|\le \sum_{\bk \ge \bzero} A(\bk+\bun)^\bn|\bz^\bk|&\le& A\prod_{1\le m\le g}\,\,\sum_{k_m\ge 0}(k_m+1)^{n_m}|z_m|^{k_m}\\
&\le& \frac{\bn !A}{\prod_m(1-|z_m|)^{n_m+1}}\\
&=&\frac{\bn!A}{(\bun-\bB \bz  \bB)^{\bn +\bun}}\le  \frac{\bn !A}{(1-|\bz|_\infty)^{g+|n|_1}}.
\end{eqnarray*}

Soit $\bk =(k_1,\ldots, k_g)$. Posant pour tout $m$ 

$$u_m=\frac{k_m+1+(n_m+1)|z_m|}{k_m+n_m+2},$$
\noindent l'int{\'e}grale de Cauchy donne

\begin{eqnarray}\label{equation:plusCauchy}
|f^{(\bk)}(\bz)|&=&\left|\frac{\bk !}{(2\pi i)^g}\int_{|\zeta_1|=u_1}\int_{|\zeta_2|=u_2}\dots
\int_{|\zeta_g|=u_g} \frac{f(\zeta_1, \zeta_2, \dots, \zeta_g)}{\prod_m (\zeta_m-z_m)^{k_m+1}}
d\zeta_1d\zeta_2\dots d\zeta_g  \right|\nonumber \\
&\le & A\bn!\bk!\frac{(\bk+\bn+\bdeux)^{\bk+\bn+\bdeux}}{(\bk+\bun)^{\bk+\bun}(\bn+\bun)^{\bn+\bun}}
\frac{1}{\prod_m(1-|z_m|)^{n_m+k_m+2}}\nonumber\\
&=& A\bn!\bk!\frac{(\bk+\bn+\bdeux)^{\bk+\bn+\bdeux}}{(\bk+\bun)^{\bk+\bun}(\bn+\bun)^{\bn+\bun}(\bun -\bB \bz \bB)^{\bn +\bk+\bdeux}}
\nonumber\\
&\le & A\bn!\bk!\frac{(\bk+\bn+\bdeux)^{\bk+\bn+\bdeux}}{(\bk+\bun)^{\bk+\bun}(\bn+\bun)^{\bn+\bun}(1-|\bz|_\infty)^{2g+ |n|_1+ |k|_1}}
\end{eqnarray}

\'Etant  donn{\'e} $\bc = (c_1, \ldots, c_m)\in P(\bzero,1)$ il est alors naturel de consid{\'e}rer la s{\'e}rie {\it recentr{\'e}e }
en $\bc$ \index{s{\'e}rie recentr{\'e}e}

$$F_\bc (\by)=F_\bc (y_1, \ldots, y_g)=f(\bc+\by\star(\bun -\bB \bc \bB))=f((c_m +y_m (1-|c_m|))_m)$$
\noindent d{\'e}finie pour $\by \in P(\bzero, 1).$

Reprenons l'in{\'e}galit{\'e} (\ref{equation:plusCauchy}) et d{\'e}duisons

\begin{lemme}[Recentrage]\label{lemme:plusrecentrage}
Soit $g\ge 1$ un entier, $A\ge 1$ un r{\'e}el et $\bn \ge \bun$ dans
$\NN^g$ et soit $f=\sum_{\bk\ge \bzero}
f_\bk\bz^\bk$ une s{\'e}rie enti{\`e}re d'ordre de grandeur $(A, \bn)$. Soit $\bc \in P(\bzero,1)$
et notons $F_\bc (\by)=f(\bc+\by\star(\bun -\bB \bc \bB))$ la s{\'e}rie recentr{\'e}e
de $f$ en $\bc$. Alors $F_\bc$ est une s{\'e}rie enti{\`e}re d'ordre de grandeur $(A_\bc,  \bn+\bun)$
avec

$$A_\bc = \bn !A\exp (g+|\bn|_1)2^{g+|\bn|_1}(\bun-\bB \bc \bB)^{-\bn-\bdeux}$$
\end{lemme}

Pour tout entier positif  $u$  on note $R_u(\bz)$\index{$R_u(\bz)$, le reste d'ordre $u$ d'une  s{\'e}rie} le reste de la s{\'e}rie

$$f(\bz)=\sum_{|\bk|_1 \le u-1 }f_\bk\bz^\bk + R_u(\bz).$$

Soit $0<R<1$ et $\bz \in P(\bzero, R)$  on a

\begin{eqnarray*}
|R_u(\bz)|&=& \left| \sum_{|\bk|_1\ge u}  f_\bk\bz^\bk\right| \le   \sum_{\stackrel{|\bk|_1\ge u}{\bk\le (u-1)\bun}}  \left|
f_\bk\bz^\bk\right| +
 \sum_{\stackrel{|\bk|_1\ge u}{\bk\not\le (u-1)\bun}}   \left| f_\bk\bz^\bk\right|  \\
&\le&\sc u^g|\bz|_\infty^{u}Au^
{|\bn|_1}\,\,+\,A\sum_{1\le m\le g}\,\,\,\,\sum_{k_m\ge u} (k_m+1)^{n_m}|z_m^{k_m}|
\times \prod_{\stackrel{1\le t \le g}{t\not = m}}\sum_{k_t\ge 0} (k_t+1)^{n_t}|z_t^{k_t}|\\
&\le &u^g|\bz|_\infty^{u}Au^
{|\bn|_1}+\frac{\bn !A}{\prod_m(1-|z_m|)^{n_m+1}}\sum_{1\le m\le g} |z_m|^u(1+u)^{n_m}\\
&\le&u^g|\bz|_\infty^{u}Au^{|\bn|_1}+\frac{\bn!A}{(\bun-\bB \bz  \bB)^{\bn +\bun}}\times g\times |\bz|_\infty^u(1+u)^{|\bn|_\infty}  \\
&\le&  u^g|\bz|_\infty^{u}Au^{|\bn|_1}+\frac{\bn !gA|\bz|_\infty^u(1+u)^{|\bn|_\infty}}{(1-|\bz|_\infty)^{g+|\bn |_1}}\\
&\le& u^g|\bz|_\infty^{u}Au^{ng}+\frac{1}{2} B(1+u)^{n}|\bz|_\infty^u\le B(u+1)^{(n+1)g}|\bz|_\infty^u
\end{eqnarray*}
\noindent avec $n=|\bn|_\infty$ et $B=\frac{\bn!2 Ag}{(1-R)^{g+|\bn|_1}}$.

\begin{lemme}[Majoration du reste]\label{lemme:majorationdureste}
Soit $f(\bz)$ une s{\'e}rie de $g$ variables d'ordre de grandeur $(A, \bn)$. Soit $u$ un entier positif ou nul et
$R_u(\bz)$ le reste d'ordre $u$. Soit $R$ un r{\'e}el
strictement compris entre $0$ et $1$. Pour tout $\bz$ dans $P(\bzero, R)$ on a

$$|R_u(\bz)|\le B(u+1)^{(n+1)g}|\bz|_\infty^u$$
\noindent avec $n=|\bn|_\infty$ et $B=\frac{\bn! 2Ag}{(1-R)^{g+|\bn|_1}}$.

En outre,  si $0<\kappa<1$ est un r{\'e}el et si 

$$u\ge \max (\frac{16(ng)^2}{(\log R)^2},\frac{2(\log \kappa -\log B)}{\log R})$$
\noindent alors $|R_u(\bz)|\le \kappa$ pour $\bz \in P(\bzero, R)$. 

\end{lemme}
 
En effet, si $u\ge \frac{16(ng)^2}{(\log R)^2}$ alors 
$(n+1)g\sqrt u \le \frac{u|\log R |}{2}$ donc $\log B + (n+1)g\log (1+u) +u\log R
\le \log B +\frac{u\log R}{2}\le \log \kappa$ car   $u \ge \frac{2(\log \kappa -\log B)}{\log R}$. \hfill $\Box$

\subsection{Prolongement analytique sur un disque}\label{subsection:prol}

Dans ce paragraphe  on veut montrer qu'une s{\'e}rie enti{\`e}re $f$ d'une variable et 
d'ordre de grandeur $(A,n)$ avec $A\ge 1$ et $n \ge 1$,
major{\'e}e en module par un  $\epsilon >0$
sur un petit disque $D=D(c,r)$ inclus dans $D(0,1)$
peut {\^e}tre agr{\'e}ablement major{\'e}e en module  sur le gros
disque $D(0,\frac{1}{2})$. 

On proc{\`e}de par prolongement {\`a} des disques de plus en plus gros.

On introduit donc la 

\begin{definition}[Fils d'un disque équilibré]
Un disque ouvert non vide  contenu
dans le disque unit{\'e} est dit {\rm {\'e}quilibr{\'e}} si et seulement si sa distance au
cercle  unit{\'e} est {\'e}gale {\`a} son rayon.
Pour $c\in D(0,1)$, on note  $D_c$ le disque {\'e}quilibr{\'e} de centre
$c$.\index{$D_c$, le disque {\'e}quilibr{\'e} de centre
$c$} Sont rayon est  $r=\frac{1-|c|}{2}$.

Si de plus 
$|c|>\frac{1}{5}$, posons $c'=(|c|-\frac{r}{2})\times \frac{c}{|c|}$ et soit $D_{c'}$ le disque {\'e}quilibr{\'e} de centre
$c'$. Alors $|c|- \frac{r}{2} = \frac{5|c|-1}{4}>0$
donc le rayon $r'$ de $D_{c'}$ v{\'e}rifie $r'=\frac{5}{4}r$. De plus
$1-|c'|=\frac{5}{4}(1-|c|)$.
On dit que $D_{c'}$ est le fils\index{fils d'un disque} de $D_c$. Si $|c|\le \frac{1}{5}$ alors
le fils de $D_c$ est par d{\'e}finition  $D_0=D(0,\frac{1}{2})$.
\end{definition}

\begin{figure}[htb]
\begin{center}
\epsfig{file=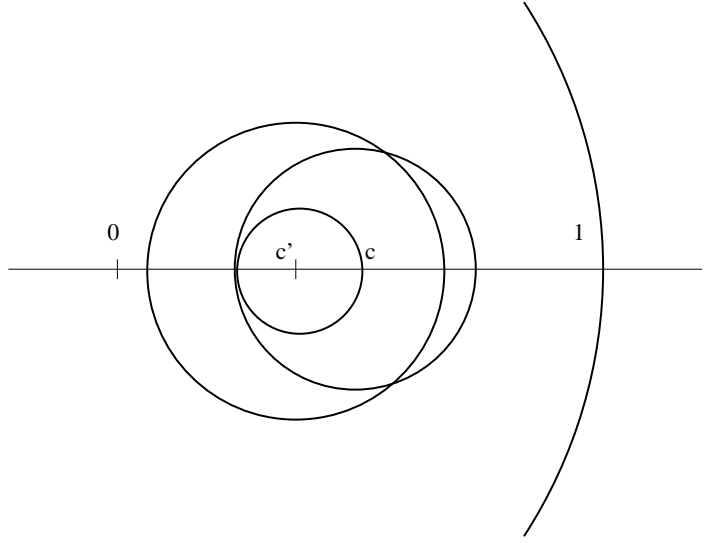}
\end{center}
\caption{Fils d'un disque {\'e}quilibr{\'e}}\label{figure:cercles}
\end{figure}

Soit $f(z)$ une s{\'e}rie enti{\`e}re d'une variable et d'ordre de grandeur $(A,n)$ avec $A\ge 1$ et $n\ge 1$.
Soit $D_c\subset D(0, 1)$ un disque {\'e}quilibr{\'e} o{\`u} $f$ est major{\'e}e en module par $0<\epsilon <2^{-100}$. 
On note $r$ le rayon de $D_c$.
Soit $D_{c'}$ le disque fils de $D_c$ et $r'$ son rayon. Le disque  de centre $c'$ et de rayon
$r/2$, est contenu dans $D_c$. Donc $f$ y est major{\'e}e en module par
$\epsilon$.

La formule de Cauchy majore  les d{\'e}riv{\'e}es de $f$ en $c'$.

\begin{equation*}
|f^{(k)}(c')|=\left|\frac{k!}{(2\pi i)}\int_{|\zeta|=r/2}
\frac{f(c' + \zeta)}{\zeta^{k+1}}d\zeta \right|\le \epsilon\frac{2^{k}k!}{r^{k}}
\end{equation*}

On veut majorer $|f|$ sur $D_{c'}$. On choisit un entier positif $u$ et on majore s{\'e}par{\'e}ment
la partie principale d'ordre $u$ en $c'$ not{\'e}e $P_{c',u}$ et le reste $R_{c',u}=f-P_{c',u}$.

D'une part

\begin{eqnarray}\label{equation:majorprincipal}
|P_{c',u}(c'+z)| &\le & \sum_{0\le k \le u-1} \left| \frac{f^{(k)}(c')}{k !}z^k\right|
\nonumber\\ &\le &   \sum_{0\le k \le u-1} \left(\frac{5r}{4}\right)^{k}\frac{2^{k}\epsilon}{r^{k}}\le \epsilon \left(\frac{5}{2}\right)^{u}.
\end{eqnarray}

D'autre part, posant $z=c'+ y (1 -|c'|)$, le reste $R_{c',u}(c' + y (1 -|c'|))$ n'est autre 
que le reste d'ordre $u$ en $0$ de la s{\'e}rie recentr{\'e}e $y \mapsto F_{c'}(y)$.
Puisque $z=c'+ y (1 -|c'|)$ appartient au disque {\'e}quilibr{\'e} $D_{c'}$,
le vecteur $y$ parcourt le disque  {\'e}quilibr{\'e} $D_0$. En d'autres termes, 
$|y| \le \frac{1}{2}$. On applique les lemmes \ref{lemme:plusrecentrage} de recentrage
et \ref{lemme:majorationdureste} de majoration du reste. 

La s{\'e}rie recentr{\'e}e $F_{c'}$ est d'ordre de grandeur $(A_{c'}, n+1)$ avec
$A_{c'}$  major{\'e}e par $n! A\left( \frac{2e}{1-|c'| }\right)^{n+2}\le n!A\left( \frac{e}{r} \right)^{n+2}$.
 
Suivant les notations du lemme \ref{lemme:majorationdureste} on pose 

$$B_{c'}=\frac{2(2e)^{n+2}(n+1)!n!A}{r^{n+2}}
\le A\exp( \cO n^2(1+|\log r|)).$$

Pour $z \in D_{c'}$ on a 

$$|R_{c',u}(z)|\le B_{c'} (1+u)^{n+2}2^{-u}$$

Soit alors $\kappa >0$ le r{\'e}el tel que $\log \kappa =\frac{\log 2\log \epsilon}{12\log \frac{5}{2}}$
et soit $u$ le plus petit entier plus grand que $\frac{4|\log \kappa|}{\log   2}= \frac{ |\log \epsilon |}{3\log \frac {5}{2}}$.

 On suppose que

$$|\log \kappa |\ge |\log B_{c'}|$$
\noindent 
 donc $u\ge \frac{2(|\log \kappa| +|\log B_{c'}|)}{\log 2}$. On suppose en outre
que 

$$u\ge \frac{16(n+1)^2}{(\log 2)^2}.$$

 Alors 

$$|R_{c', u}|\le \kappa.$$

 On  d{\'e}duit de (\ref{equation:majorprincipal}) que 

$$\log |P_{c',u}(c'+z)|  \le u\log\frac{5}{2} +\log \epsilon .$$

On a  $u\le 1+\frac{ |\log \epsilon |}{3\log \frac {5}{2}}$  donc
$u\log \frac{5}{2} \le \log \frac{5}{2}+\frac{|\log \epsilon |}{3} \le \frac{|\log \epsilon |}{2}$ car 
 $|\log \epsilon | \ge 100 \log 2$.
Donc $\log |P_{c',u}(c'+z)|  \le \frac{\log \epsilon }{2}$.

Ainsi $\log |f| \le \log \left( 2\max (|P_{c',u}|,  |R_{c',u}|) \right) \le  \frac{\log 2}{12(\log 5-\log 2)}\log \epsilon + \log 2
\le 0.05\log \epsilon$ car $\log \epsilon \le -100\log 2$.

\begin{lemme}[Prolongement  au disque  fils]\label{lemme:prolongementaufils}
Il existe une constante positive effective $c_{12}$ telle que l'{\'e}nonc{\'e} suivant soit vrai :

Soit  $f(z)$ est  une s{\'e}rie enti{\`e}re d'une  variable et d'ordre de grandeur 
$(A, n)$ avec $A\ge 1$ et $n \ge 1$. Soit $D=D_c\subset
D(0, 1)$ un disque {\'e}quilibr{\'e} de centre $c$ et de rayon $r$
o{\`u} $f$ est major{\'e}e en module par un $0 <\epsilon < 1$. 
Soit $D_{c'}$ le disque fils de $D_c$. On suppose que  $-\log \epsilon$ est minor{\'e} 
 par $c_{12}(\log A + n^2|\log r|)$. 
    Alors $f$ est major{\'e}e en module
sur le disque fils  $D_{\bc'}$  par $\epsilon ^\frac{1}{20}$.
\end{lemme}

\subsection{Polygone de Newton d'une s{\'e}rie enti{\`e}re}\label{subsection:newton}

Soit $f$ une s{\'e}rie enti{\`e}re d'une variable complexe
 et soit
$R\le \infty $ son rayon de convergence, suppos{\'e} non nul.
Soit $r$ un r{\'e}el positif inf{\'e}rieur {\`a} $R$.
 On note $D=D(0,r)$.

On s'int{\'e}resse
aux  z{\'e}ros de
$f$ dans $D$. Combien sont ils ? O{\`u} sont ils ? Comment sont ils affect{\'e}s par une 
petite perturbation de $f$ ?

Pour r{\'e}pondre {\`a} ces questions, on cherche {\`a} enfermer les z{\'e}ros de $f$ dans une collection
finie de petits disques disjoints tels que le nombre de z{\'e}ros de $f$ dans chaque
disque ne soit pas modifi{\'e} par une petite perturbation.

On {\'e}tudie d'abord la situation autour de z{\'e}ro. On suppose que $f(0)=1$ donc
$f(z)=1+\sum_{k \ge 1} f_k z^k$. On note $d$ le degr{\'e} de $f$ en $z$, qui est en g{\'e}n{\'e}ral
infini. Le {\it nuage de Newton} associ{\'e} {\`a} $f$ est l'ensemble de points $(k,-\log |f_k|)$ pour 
$ k \ge 0$ et $f_k \not =0$. Le {\it polygone de Newton }\index{polygone de Newton} de $f$ est la fonction $\cN$
 de $[0,d]$
dans $\RR$ d{\'e}finie comme le maximum des fonctions affines
  $\phi : [0,d]\rightarrow \RR$
  qui passent en dessous du nuage de Newton (c'est-\`a-dire
  $\phi(k)\le -\log |f_k|$ pour tout $k$). La fonction $\cN$  est bien d{\'e}finie car le rayon de convergence 
$R$ est non nul. C'est une fonction  convexe de $[0, d]$.
Elle est affine sur tout intervalle ouvert d{\'e}limit{\'e} par deux entiers cons{\'e}cutifs.

En effet, soit $k\ge 0$ un entier inf{\'e}rieur {\`a} $d$. Pour tout $\epsilon >0$ il existe une fonction affine
$\phi$
qui passe sous le nuage de Newton  et 
 telle que $\cN(k)-\epsilon\le \phi (k) \le \cN(k)$. De m{\^e}me 
il existe une fonction affine
$\psi$ qui passe sous le nuage de Newton 
 et  telle que $\cN(k+1)-\epsilon \le \psi (k+1) \le \cN(k+1)$. On d{\'e}finit la fonction
affine $\kappa_\epsilon$ de la fa{\c c}on suivante. Si $\phi(k+1) < \psi(k+1)$ et
$\psi(k) < \phi(k)$ alors $\kappa_\epsilon$ est la fonction affine qui vaut 
$\phi(k)$ en $k$ et $\psi(k+1)$ en $k+1$. Si $\phi(k+1) \ge  \psi(k+1)$ alors 
$\kappa _\epsilon = \phi$. Si $\phi(k+1)<  \psi(k+1)$ et $\psi(k)\ge \phi(k)$
alors $\kappa _\epsilon = \psi$. On v{\'e}rifie que $\kappa_\epsilon$
passe
sous  le nuage 
de Newton. Quand  $\epsilon$ tend vers $0$  la famille des $\kappa_\epsilon$ 
 converge simplement sur le segment $[k,k+1]$  vers la fonction affine $\kappa_0$
qui vaut $\cN(k)$ en $k$ et $\cN(k+1)$ en $k+1$. Donc $\kappa_0$ minore $\cN$ sur cet
intervalle. Un argument de convexit{\'e} montre qu'on a {\'e}galit{\'e}.

Ainsi  $\cN$ est continue et affine par morceaux sur $[0,d]$.

 Les {\it sommets} du polygone de Newton sont les discontinuit{\'e}s de $\cN'$ plus $(0,0)$
et {\'e}ventuellement $(d, \cN(d))$. 

Soit $k$ un entier entre $0$ et $d$. On pose $l=\cN(k)$. On appelle
tangente en $(k, l)$ au polygone de Newton, toute droite passant par $P = (k, l)$ 
et qui passe sous  le nuage de Newton. On note $\alpha^-$ la d{\'e}riv{\'e}e {\`a} gauche, qui est
la  pente de la tangente {\`a} gauche. Le vecteur $(-\alpha^-, 1)$ est orthogonal {\`a} cette droite
et tourn{\'e} vers l'int{\'e}rieur de $\cN$. De m{\^e}me $\alpha^+$ est la d{\'e}riv{\'e}e {\`a} droite.
On suppose que $\alpha^+ > \alpha^-$ donc $P$ est un sommet.
Soit $\alpha$ dans $]\alpha^-, \alpha^+[$. La position relative du nuage de  Newton et 
de la tangente en $P$ de pente $\alpha$ nous renseigne sur l'ordre de grandeur de $f(z)$ 
pour un $z$ tel que  $\log |z|=\alpha$. En effet pour tout tel $z$ et
pour tout entier positif 
$m$ on a $-\log |f_mz^m|=-\log |f_m|-m\alpha = (-\alpha,1)\cdot(m,-\log |f_m|)
\ge (-\alpha,1)\cdot (k,-\log |f_k|)$. De sorte que $f_kz^k$ est le terme dominant sur le cercle $|z|=\exp(\alpha)$.
Il reste {\`a} voir jusqu'{\`a} quel point. On se doute que si le sommet $P$ est assez anguleux, les autres termes
peuvent \^etre  n{\'e}gligeables.

Soit donc $m\not = k$ un entier positif ou nul. Le point $(m, -\log
|f_m|)$ est au dessus du polygone de Newton.
Si $m>k$ il est donc au dessus de
  la droite passant par $P$ et de pente $\alpha^+$. Donc $|f_m|\le |f_k|\exp(-(m-k)\alpha^+)$.
Donc, pour $\log |z|=\alpha$,  le terme $f_mz^m$ est major{\'e} en module par $|f_k||z|^k$ fois $\exp((m-k)(\alpha-\alpha^+))$.
La somme $\sum_{m>k} |f_k||z|^k$ est donc major{\'e}e par $|f_k||z|^k$ fois $\frac{x}{1-x}$ en posant
$x=\exp(\alpha-\alpha^+)$. 
Si $m<k$ le point $(m, -\log |f_m|)$ est au dessus de la droite passant par $P$ et de pente 
$\alpha^-$. Donc $|f_m|\le |f_k|\exp(-(m-k)\alpha^-)$.
Donc, pour $\log |z|=\alpha$,  le terme $f_mz^m$ est major{\'e} en module par $|f_k||z|^k$ fois $\exp((m-k)(\alpha-\alpha^-))$.
La somme $\sum_{m<k} |f_k||z|^k$ est donc major{\'e}e par $|f_k||z|^k$ fois $\frac{y}{1-y}$ 
 en posant
$y=\exp(\alpha^--\alpha)$. 

Pour toute fonction $h$  holomorphe  sur un voisinage du  disque ferm{\'e} de centre $0$ 
et de  rayon $\exp(\alpha)$, et non nulle sur le bord
 de ce disque, le nombre de z{\'e}ros de $h$ dans l'int{\'e}rieur de ce  disque est l'indice en $0$ du lacet image du
bord $\{ z,\,   |z|=\exp(\alpha)\}$ par $h$. C'est un invariant homotopique discret. 
Donc dans toute famille continue de fonctions
holomorphes sur un voisinage de ce  disque et jamais nulles sur son bord, le nombre de z{\'e}ros dans l'int{\'e}rieur du disque
 est
constant.

On en d{\'e}duit que si $\frac{x}{1-x}+\frac{y}{1-y}<1$  le nombre de z{\'e}ros de $f$ dans le disque ouvert
$D(0,\exp(\alpha))=\{z,\,  |z| < \exp(\alpha)\}$ est le m{\^e}me que celui de $z^k$ soit $k$ z{\'e}ros.
Cette condition est satisfaite si $x$ et $y$ sont inf{\'e}rieurs {\`a} $\frac{1}{3}$.

On appelle {\it pente } du polygone de Newton,  une valeur de la  d{\'e}riv{\'e}e de $\cN$ en un point ou elle
est d{\'e}rivable.

\begin{lemme}[Polygone de Newton]\label{lemme:Newton}
Soit $f=1+\sum_{k\ge 1}f_kz^k$ une s{\'e}rie enti{\`e}re de rayon de convergence $R>0$. Soit $\xi \in D(0,R)$ un z{\'e}ro de
$f$. Il existe une pente $\sigma$ du polygone de Newton telle que $|\log |\xi| - \sigma| \le \log 3$.

On note $\cP_3$ l'intervale $]-\infty, \log R[$ priv{\'e} des intervales $[\sigma-\log 3, \sigma +\log 3]$
o{\`u} $\sigma$ parcourt l'ensemble des pentes du polygone de Newton.

Si $\alpha$ est un r{\'e}el de $\cP_3$ il existe un unique sommet $P=(k, \cN(k))$ admettant une
tangente de pente $\alpha$. La fonction $f$ a exactement $k$ z{\'e}ros dans le disque ouvert
$D(0,\exp(\alpha))$.

On note $\cP_4$ l'intervale $]-\infty, \log R[$ priv{\'e} des intervales $]\sigma-\log 4, \sigma +\log 4[$
o{\`u} $\sigma$ parcourt l'ensemble des pentes du polygone de Newton.

Si $\alpha$ est un r{\'e}el de $\cP_4$  il existe un unique sommet $P=(k, \cN(k))$ admettant une
tangente de pente $\alpha$. Pour $\log |z| = \alpha$ on a $|f(z)|\ge \frac{|z|^k}{3}= \frac{\exp ( k\alpha )}{3}$.

\end{lemme}

Comme on pouvait s'y attendre, ce lemme est moins pr\'ecis que son
pendant non-archim\'edien. On ne peut pas l'utiliser directement si
les pentes sont trop proches les unes des autres. Dans ce cas, on
pourra former (par exemple) la s\'erie $g(z)=f(\sqrt z)f(-\sqrt z)$
dont les z\'eros sont les carr\'es des z\'eros de $f$. Le passage de $f$ \`a
$g$ clarifie la situation dans le voisinage du cercle unit\'e. On peut
r\'eit\'erer l'op\'eration si n\'ecessaire.

\subsection{Le plus petit z{\'e}ro d'une s{\'e}rie enti{\`e}re}\label{subsection:petitzero}

Soit $F=F_0+\sum_{k\ge 1}F_kz^k$ une s{\'e}rie non  constante de rayon au moins $1$ telle que $F_0\not =0$.
La s{\'e}rie normalis{\'e}e  $f=F/F_0$  admet au moins une pente. Soit
$\sigma_1$ la plus petite  des pentes. On suppose d'abord que
$\sigma_1$ est négatif. Soit  alors $\log r $ la borne inf{\'e}rieure  de $]\sigma_1, 0  [  \cap  \cP_3$. Si ce dernier
ensemble est vide on pose $r=1$. Si $r<1$, alors $f$ admet un z{\'e}ro de
module $\le r$. On veut montrer que si $F_0$ est petit alors $r$ est
petit ou bien $F$
est uniform\'ement petite. On suppose que $|F_0|<1$.
 Le segment $]\sigma_1, \log r [$ est couvert par des intervales  ferm{\'e}s de rayon $\log 3$
 centr{\'e}s en les  pentes du polygone de Newton. On note $\sigma_1 < \sigma_2 < \cdots$ les pentes
successives. On a $\sigma_2 \le \sigma_1+2\log 3$, $\sigma_3 \le \sigma_1+4\log 3$, \ldots,
$\sigma_k \le \sigma_1+2(k -1)\log 3$, tant que  $\sigma_1+(2k -3)\log 3<\log r$.
On pose donc 

$$\ell=\lceil \frac{\log r-\sigma_1+\log 3}{2\log 3}  \rceil$$
\noindent  et pour tout $1\le k\le \ell$
on a $\sigma_k \le \sigma_1 + 2(k-1)\log 3$ et donc $\cN(k)\le k\sigma_1+k(k-1)\log 3$. 
Cela prouve en particulier que le degr{\'e} $d$ de $F$ est au moins {\'e}gal 
{\`a} $\ell$.

On pose  $k=\ell$ et on obtient 

\begin{equation}\label{equation:minor}
\cN(\ell) \le \ell\sigma_1+\ell(\ell-1)\log 3.
\end{equation}

Le principe de la démonstration  est le suivant : on suppose $F_0$ petit.

Si $\sigma_1$ est grand   alors la s{\'e}rie $F$ est petite
car ses premiers coefficients sont petits.

Si $\sigma_1$ est petit et $\ell$ petit alors $r$ est petit : le polygone de Newton  est anguleux pr{\`e}s de l'origine
et il y  a une petite racine. 

Si $\sigma_1$ est petit  et $\ell$ grand alors la pente du polygone de Newton varie peu au d{\'e}but, et la s{\'e}rie
$F$ a de grand coefficients.

Pour formaliser ce raisonnement, nous supposons maintenant que $F$ est d'ordre de grandeur $(A, n)$ avec $A\ge 1$
et $n\ge 1$. Donc $-\log |f_k|=-\log |F_k| +\log |F_0|$
est minor{\'e}e par $-\log A -n\log (k+1) + \log | F_0 |$ qui est une fonction convexe de $k$ et qui minore
donc le polygone de Newton.  
Pour $k=\ell$ on obtient 

\begin{equation}\label{equation:major}
\cN(\ell)\ge -\log A -n\log (\ell+1) + \log |F_0|.
\end{equation}

Comme $\ell \le \frac{\log r-\sigma_1+3\log 3}{2\log 3}$ on a 
$\sigma_1 \le -2\ell\log 3+\log r +3\log 3$.
En reportant dans l'in{\'e}quation (\ref{equation:minor}) on a 
$\cN(\ell) \le -\ell ^2\log 3 +\ell (\log r +2 \log 3)$. L'in{\'e}quation (\ref{equation:major}) donne alors
$\ell ^2\log 3 -\ell (\log r +2 \log 3) \le \log A +n\log (\ell+1) - \log |F_0| \le \log A +n\ell - \log |F_0|$.
Donc $\ell$ satisfait l'in{\'e}galit{\'e} quadratique

\begin{equation}\label{equation:quad}
\ell^2\log 3-\ell (\log r+2\log 3+n)-\log A+\log |F_0| \le 0.
\end{equation}

Si $\ell \ge \log r +2\log 3+n$ alors  on d{\'e}duit de l'in{\'e}quation (\ref{equation:quad})
que 
$\ell^2(\log 3-1)\le \log A -\log |F_0|$. Au total

\begin{equation}\label{equation:majell}
\ell \le \max (\log r+2\log 3+n, \, \sqrt{\frac{\log A - \log |F_0|}{\log 3 -1}}).
\end{equation}

Comme $\ell \ge \frac{\log r-\sigma_1+\log 3}{2\log 3}$  
on a $\sigma_1 \ge -2\ell \log 3+\log r+\log 3$.
On d{\'e}duit de l'in{\'e}quation (\ref{equation:majell})

\begin{equation*}
\scriptstyle \sigma_1 \ge \min ((1-2\log 3)\log r-2n\log 3-(2\log 3)^2+\log 3, \, -2\log 3\sqrt{\frac{\log A - \log |F_0|}{\log 3 -1}}+\log r +\log 3).
\end{equation*}
\noindent On rappelle que  $r \le 1$, et on suppose que 

$$-\log |F_0| \ge \cO (\log A+ n^2).$$

 On en d{\'e}duit alors
que 

\begin{equation*}
\sigma_1 \ge -10 \sqrt{- \log |F_0|} +\log r  
\end{equation*}

Si $-\log r \ge \sqrt{- \log |F_0|}$ on s'estime heureux puisqu'on a montr{\'e} que $F$ admet un 
z{\'e}ro tr{\`e}s petit. Sinon on a $\sigma_1 \ge - 11 \sqrt{- \log
  |F_0|}$. 
On observe que cette dernière inégalité est vraie aussi si $\sigma_1$
est positif ou nul.

Donc pour tout entier $k\ge 0$
on a $f_k=\frac{F_k}{F_0}\le \exp(-k\sigma_1)\le \exp(11 k\sqrt{-\log |F_0|})$. Si $z\in D(0,\frac{1}{2})$
est un complexe de module inf{\'e}rieur {\`a} $\frac{1}{2}$ alors $|F_k||z|^k\le |F_0|\exp(k(11 \sqrt{-\log |F_0|}-\log 2))$
et pour tout entier positif $u$, la partie principale $P_u(z)=\sum_{0\le k < u}F_uz^u$ est major{\'e}e en module
par $|F_0|\exp(u(11 \sqrt{-\log |F_0|}-\log 2))$. 

Si on choisit $u=\lfloor \frac{\sqrt{-\log |F_0|}}{22}\rfloor $ alors
$\log |P_u(z)| \le \log |F_0| +u(-\log 2 +11\sqrt{-\log |F_0|})\le   \frac{\log |F_0|}{2}$.

Pendant ce temps l{\`a},  on peut majorer le reste $R_u(z)$ {\`a} l'aide du lemme \ref{lemme:majorationdureste}.
On demande que $u\ge \frac{16n^2}{(\log 2)^2}$ ce qui est acquis si 

$$-\log |F_0| \ge \cO n^4.$$

Suivant les notations du lemme \ref{lemme:majorationdureste}
on pose $B=2^{n+2}n!A$. On peut majorer $R_u(z)$ en module par $\kappa$ pourvu que
 $-\log \kappa \le  -\log B +u\frac{\log 2}{2}$.  On note que $u=\lfloor \frac{\sqrt{-\log |F_0|}}{22}\rfloor$
est plus grand que $\frac{\sqrt{-\log |F_0|}}{23}$ si  $\sqrt{-\log |F_0|} \ge \cO.$

Si $\sqrt{-\log |F_0|} \ge \cO (\log A+ n^2)$ alors 

$$\log B\le \frac{\log 2}{92}\sqrt{-\log |F_0|}$$ 
\noindent donc  $-\log B +u\frac{\log 2}{2} \ge \sqrt{-\log |F_0|} \times \frac{\log 2}{92}$.  

On pose
donc $-\log \kappa = 0.007\sqrt{-\log |F_0|}$.

Comme la partie principale $P_u(z)$ est major{\'e}e par $\sqrt{|F_0|}$ qui est plus petit que $\kappa$ on 
a $\log |F(z)|\le \log 2-0.007 \sqrt{-\log |F_0|} \le -0.006  \sqrt{-\log |F_0|}$ si  $\sqrt{-\log |F_0|} \ge \cO.$

\begin{lemme}[Plus petit z\'ero]\label{lemme:proximite}
Il existe une constante positive effective $c_{13}$ telle que l'{\'e}nonc{\'e} suivant soit vrai :

Soit $F(z)=F_0+F_1z+\cdots$ une s{\'e}rie enti{\`e}re d'ordre de grandeur  $(A,n)$ avec $A\ge 1$ et $n\ge 1$. 
 On suppose que $|F_0|< 1$ et 
$\sqrt{-\log |F_0|}\ge c_{13} (n^2+\log A)$.
Alors ou bien $f$ a un z{\'e}ro $\xi$ tel que $\log |\xi| \le -\sqrt{-\log |F_0|}$, ou bien
$f(z)$ est major{\'e}e en module pour $z\in D(0,\frac{1}{2})$ par $\kappa $ tel que
$\log \kappa =\frac{-\sqrt{-\log |F_0|}}{200}$.
\end{lemme}

\subsection{Stabilit{\'e} des z{\'e}ros d'une s{\'e}rie enti{\`e}re}\label{subsection:stabilitezeros}

On se donne maintenant une s{\'e}rie enti{\`e}re $f$  d'ordre de grandeur 
$(A, n)$ avec $A\ge 1$ et $n\ge 1$ et un complexe $c$ tel que $|c|<1$.
On veut montrer que si $f(c)$ est petite alors $c$ est proche d'un z{\'e}ro
de $f$ ou bien $f$ est petite sur le disque $D(0,\frac{1}{2})$.

On note $F_c(y)=f(c+y(1-|c|))=F_0+F_1y+\cdots$ la s{\'e}rie recentr{\'e}e en $c$.
Elle est d'ordre de grandeur $(A_c, n+1)$ avec $A_c=An!\exp(n+1)2^{n+1}
(1-|c|)^{-n-2}$. Appliquons le lemme \ref{lemme:proximite} {\`a} la s{\'e}rie 
recentr{\'e}e $F_c$. On suppose donc que $F_0=f(c)$ v{\'e}rifie $|f(c)|<1$ et

$$\sqrt{-\log |f(c)|}\ge c_{13}(n^2  + \log A_c).$$

Il vient que $f$ a un z{\'e}ro $\xi$ tel que $\log |c-\xi|\le -\sqrt{-\log |f(c)|}$
ou bien $f$ est major{\'e}e en module par $\kappa$ sur 
le disque {\'e}quilibr{\'e} $D_c$ o{\`u}  $\log \kappa
=-\frac{\sqrt{-\log |f(c)|}}{200}$. Dans ce dernier cas, on peut appliquer
le lemme \ref{lemme:prolongementaufils}.  Soit
$w=\lceil \frac{-log(1-|c|)}{\log \frac{5}{4}} \rceil$ et
soit $\nu$ tel que $\log \nu =\frac{\log \kappa}{20^w}$
et supposons que
  
$$-\log \nu  \ge  c_{12}(\log A + n^2|\log \left(\frac{1-|c|}{2}\right)|).$$ 

On applique $w$ fois le lemme \ref{lemme:prolongementaufils}
et on montre que  $f$ est major{\'e}e en module par $\nu$ sur le disque $D(0,\frac{1}{2})$.

On note que $20^w \le 20(1-|c|)^{\frac{\log 20}{\log \frac{4}{5}}} \le 20(1-|c|)^{-14}$.

\begin{lemme}[Stabilit{\'e} d'un z{\'e}ro]\label{lemme:stabunzero}
Il existe une constante effective positive $c_{14}$ telle  que l'{\'e}nonc{\'e} suivant soit vrai :

Soit $f$ une s{\'e}rie enti{\`e}re d'ordre de grandeur $(A, n)$ avec $A\ge 1$ et $n\ge 1$. 
Soit $c$ dans $D(0,1)$. On suppose que $\epsilon = |f(c)| < 1$. Soit
$\nu $ tel que $\log \nu = -\frac{(1-|c|)^{14}\sqrt{-\log \epsilon }}{c_{14}}$.

On suppose que $-\log \nu  \ge  c_{14}(\log A + n^2(1+|\log ({1-|c|})|)).$

Alors $f$ admet un z{\'e}ro $\xi$ tel que $\log |c-\xi|\le -\sqrt{-\log |f(c)|}$
ou bien $f$ est major{\'e}e en module par $\nu$ sur 
le disque $D(0,\frac{1}{2})$.
\end{lemme}

Puisque les zéros de $f$ bougent peu sous l'effet d'une petite perturbation,
ils ne doivent pas beaucoup s'éloigner des zéros de la partie
principale. Et c'est un moyen commode de les localiser.

Soit donc $f$ une s{\'e}rie enti{\`e}re  d'ordre de grandeur $(A,n)$ 
avec $A\ge 1$ et $n\ge 1$. Soit $r$ un r{\'e}el entre $1/2$ et $1$. 
Soit $\epsilon $ un  r{\'e}el  entre $0$ et $1$ et soit $\nu$ un r{\'e}el positif tel que
$\log \nu  = -\frac{(1-r)^{14}\sqrt{-\log \epsilon}}{c_{15}}$ avec $c_{15}$ constante
plus grande que $1$ et $c_{14}$. 
On suppose que 
$|f|$  n'est pas major{\'e}e par $\nu+\epsilon$ sur $D(0,\frac{1}{2})$.
On suppose que $-\log \nu  \ge   c_{15}(\log A + n^2(1+|\log ({1-r})|))$.

On cherche 
un entier positif $u$ tel que le reste $R_u(z)$ soit major{\'e} 
en module par $\epsilon$ sur $D(0,r)$. Selon le lemme
\ref{lemme:majorationdureste} il faut que $u\ge \frac{2(\log \epsilon - \log B)}{\log r}$
avec $B = \frac{n!2A}{(1-r)^{n+1}}$. Si $c_{15}$ est assez grand alors

$$|\log B| \le |\log \epsilon|$$
\noindent
donc il suffit que $u\ge \frac{4|\log \epsilon|}{|\log r|}$. Soit donc $u$ le plus petit
entier satisfaisant cette condition. On v{\'e}rifie que $u\ge \frac{4|\log \epsilon|}{|\log r|} \ge 
\frac{16n^2}{(\log r)^2}$ si $c_{15}$ est assez grande. Donc
 $R_u(z)$ est major{\'e} par $\epsilon$ en module sur $D(0,r)$. 
On pose $\rho = \exp(-\sqrt{-\log \epsilon})$.

La partie principale $P_u(z)$ est un polyn{\^o}me de degr{\'e} $u-1$
qui a donc $u-1$ z{\'e}ros dans $\CC$. Donc l'intervalle $[r-4u\rho,r]$ 
contient\footnote{Si $c_{15}$ est assez grand alors $4u\rho$ est (beaucoup) plus petit
que $r$, donc l'intervalle en question est constitu\'e de r\'eels
positifs.} au moins un r{\'e}el positif 
$R$ tel que $|R-|\xi||>2\rho$ pour tout
z{\'e}ro $\xi$ de $P_u$.
Soit $\cD$ le ferm{\'e} de $\CC$ obtenu en retirant au disque ferm{\'e} 
$\bar D(0, R)$ tous les disques ouverts $D(\xi, 2\rho)$ o{\`u} les $\xi$
sont les z{\'e}ros de $P_u(z)$. D'apr{\`e}s le lemme \ref{lemme:stabunzero} le polyn{\^o}me $P_u(z)$ est strictement minor{\'e} en
module par $\epsilon$ sur le  domaine ferm{\'e} $\cD$. Donc $f(z)$ et 
$P_u(z)$ n'ont pas de z{\'e}ro dans $\cD$. Elles ont le m{\^e}me  nombre de z{\'e}ros 
dans $D(0,R)$. Elles ont le m{\^e}me  nombre de z{\'e}ros 
dans chaque $D(\xi, 2\rho)$.  Donc les z{\'e}ros de $P_u(z)$ dans $D(0,R)$
approchent ceux de $f(z)$ {\`a} distance  $4\rho$. On obtient le

\begin{lemme}[Stabilit{\'e} globale]\label{lemme:stabglobalezeros}
Il existe une constante effective positive $c_{16}$ telle que l'{\'e}nonc{\'e} suivant soit vrai :

Soit $f$ une s{\'e}rie enti{\`e}re d'ordre de grandeur $(A, n)$ avec $A\ge 1$ et $n\ge 1$. 
Soit $r$ et $\rho$ deux  r{\'e}els tels que $\frac{1}{2}\le r<1$ et
$0<\rho <1$. Soit $u$ le plus
petit entier plus grand que $\frac{4(\log \rho)^2}{|\log r|}$. 
On suppose que $-(1-r)^{14}\log \rho\ge c_{16}(\log A + n^2|\log
({1-r})|))$.
Alors
$f$ satisfait 
l'une au moins des deux propri\'et\'es suivantes :

\begin{enumerate}
\item Sur  le disque $D(0,\frac{1}{2})$, le logarithme 
$\log |f|$ du module de $f$  est major\'e par $\frac{(1-r)^{14}\log \rho}{c_{16}}$.
    
\item   Il existe un r{\'e}el positif $R$ tel que $r-4u\rho \le R \le r$
et tel que dans le disque $D(0,R)$ les z{\'e}ros de $f(z)$ sont approch{\'e}s {\`a} 
distance $4\rho$  par ceux de la partie principale $P_u(z)$ de degr{\'e} $u-1$. En particulier, il y a au plus 
 $u$
tels  z{\'e}ros.
\end{enumerate}
\end{lemme}

\begin{theoreme}[Z{\'e}ros d'une s{\'e}rie]
Soit $f=\sum_{r\ge 0}f_rz^r$ une s{\'e}rie enti{\`e}re d'ordre de grandeur $(A, n)$ avec $A\ge 1$ et $n\ge 1$. Soit $R$
un r{\'e}el compris strictement entre $0$ et $1$. Soit
$\mu $ la partie positive
de  $ -\log \max_{|z|\le \frac{1}{2}} |f(z)|$.
Le nombre de z{\'e}ros de module $\le R$ est polynomial en $n$, $\log A$, $(1-R)^{-1}$
et $\mu$.

Il existe un algorithme  qui pour  $f$ et $R$ comme
ci-dessus\footnote{La s\'erie $f$ est donn\'ee sous la forme d'un oracle  qui 
calcule  les   coefficients $f_r$  en temps polynomial en $r$ et en la pr{\'e}cision absolue requise.}  et
pour $k$ entier positif,  retourne 

\begin{itemize}
\item un rationnel $R'$ tel que $|R-R'|<\exp(-k)$,

\item le nombre de z\'eros de $f$ dans $D(0,R')$,

\item une approximation de ces z\'eros  \`a $\exp(-k)$ pr\`es,
\end{itemize}
\noindent  en temps d\'eterministe  
polynomial en $n$, $\log A$, $(1-R)^{-1}$, $\mu$ et la pr{\'e}cision
absolue $k$ requise. 
\end{theoreme}

Cela d{\'e}coule du lemme \ref{lemme:stabglobalezeros}.
Il suffit de rappeler l'existence de tels algorithmes pour la recherche des racines d'un polyn{\^o}me.
\hfill $\Box$

\backmatter

\bibliography{couveignes}
\begin{theindex}
\indexentry{$S^gX$, la puissance sym{\'e}trique $g$-i{\`e}me de $X$}{1}
\indexentry{$J_X$, la jacobienne de $X$}{1}
\indexentry{$\cH$, le demi-plan de Poincar{\'e}}{4}
\indexentry{$\cH^*$, le demi-plan de Poincar{\'e} avec les pointes}{4}
\indexentry{$\nu_2$, le nombre de points elliptiques d'ordre $2$}{4}
\indexentry{$\nu_3$, le nombre de points elliptiques d'ordre $3$}{4}
\indexentry{$\cP_2$, le diviseur des points elliptiques d'ordre $2$}{4}
\indexentry{$\cP_3$, le diviseur des points elliptiques d'ordre $3$}{4}
\indexentry{$q$, $q_0$, $q_\infty$, $P_0$, $P_\infty$}{4}
\indexentry{$W(P)$, l'involution d'Atkin-Lehner appliqu{\'e}e au point $P$}{4}
\indexentry{$D(x,r)$, le disque de centre $x$ et de rayon $r$}{4}
\indexentry{$R_\infty = 0.005$ }{4}
\indexentry{$R_0 = 1-\frac{1}{p}$ }{4}
\indexentry{$\Delta_d$, le diviseur $(d-1)(0) + (d-1)(\infty)+\lfloor \frac{d}{2}\rfloor \cP_2 + \lfloor \frac{2d}{3}\rfloor \cP_3$}{4}
\indexentry{$\cH^{d}(  \Delta_d)$, l'espace des formes diff{\'e}rentielles de degr{\'e} $d$ et de diviseur $\ge -\Delta_d$}{4}
\indexentry{$\cD_d$, l'ensemble des $\omega =(2i\pi)^d f(q)(d\tau)^d=\frac{f(q)}{ q^d}(dq)^d$ o{\`u} $f(q)$ est une forme modulaire parabolique primitive, propre et normalis{\'e}e  sur $\Gamma=\Gamma_0(p)$ et de poids $2d$}{5}
\indexentry{$H=\cH^1\oplus \bar \cH^1$, l'espace des  formes harmoniques}{7}
\indexentry{$D_\epsilon $, le produit $D_\epsilon =  D(0,R_{\epsilon_1})\times \cdots \times  D(0,R_{\epsilon_g})$ pour un $\epsilon = (\epsilon_k)_{1\le k \le g}\in \{0,\infty \}^g$}{9}
\indexentry{$\mu_\epsilon$ pour $\epsilon =  (\epsilon_k)_{1\le k \le g}\in \{0,\infty \}^g$}{9}
\indexentry{$\cJ_{\cD_1,\epsilon}$, le  d{\'e}terminant jacobien de l'application $\mu_\epsilon$ en $(q_{\epsilon_1,1},\ldots,q_{\epsilon_g,g})$}{9}
\indexentry{$\cS_\infty(F)$,  $\cS_0(F)$,  $|F|$, $\hat \cS_\infty(F)$,
   et $\hat \cS_0(F)$, les cinq normes de la forme $F$}{9}
\indexentry{$W_\bff $, le wronskien associ{\'e} {\`a} la famille $\bff$}{10}
\indexentry{$\cJ_\bff $, le jacobien associ{\'e} {\`a} la famille $\bff$}{10}
\indexentry{ $W_\infty (q)$, le wronskien associ{\'e} aux d{\'e}veloppements de Fourier {\`a} l'infini des formes diff{\'e}rentielles}{11}
\indexentry{$\cJ_{\infty}(\bq)$, le jacobien associ{\'e} {\`a} la famille des d{\'e}veloppements {\`a} l'infini des formes diff{\'e}rentielles}{11}
\indexentry{$g_2 = 3g-1+\nu_2+\nu_3$, la dimension de $\cH^2(\Delta_2)$}{12}
\indexentry{$W_{2,\infty}(q)$, le wronskien associ{\'e} {\`a} la famille des d{\'e}veloppements {\`a} l'infini des formes quadratiques}{12}
\indexentry{$\cJ_{2,\infty}(\bq)$, le jacobien  associ{\'e} {\`a} la famille des d{\'e}veloppements {\`a} l'infini des formes quadratiques}{12}
\indexentry{$m_2=\frac{g_2(g_2+3)}{2}$, le degr{\'e} du wronskien quadratique}{12}
\indexentry{$\cS_\infty(H)$,  $\cS_0(H)$,  $|H|$, $\hat \cS_\infty (H)$  et $\hat \cS_0(H)$,   les cinq  normes de la forme quadratique $H$}{13}
\indexentry{$\cA(\br,\chi,M)$, sous ensemble discret du tore complexe}{15}
\indexentry{r{\'e}solution d'une famille de points}{16}
\indexentry{$\star$}{20}
\indexentry{$P(\bx,\br)$, le polydisque de centre $\bx$ et de polyrayon $\br$}{20}
\indexentry{ordre de grandeur d'une s{\'e}rie}{20}
\indexentry{s{\'e}rie recentr{\'e}e}{21}
\indexentry{$R_u(\bz)$, le reste d'ordre $u$ d'une  s{\'e}rie}{21}
\indexentry{$D_c$, le disque {\'e}quilibr{\'e} de centre $c$}{22}
\indexentry{fils d'un disque}{22}
\indexentry{polygone de Newton}{24}
\end{theindex}

\end{document}